\documentclass[a4paper,11pt]{article}
\usepackage[T1]{fontenc}
\usepackage[round]{natbib}
\usepackage[table]{xcolor}
\usepackage{colortbl}
\usepackage[bookmarksnumbered,colorlinks,bookmarks,citecolor=blue,linkcolor=blue,urlcolor=blue,breaklinks,linktocpage]{hyperref}
\usepackage{amssymb}
\usepackage{amsmath,mathabx,amsthm}
\usepackage{pdflscape}
\usepackage{adjustbox}
\usepackage{afterpage}
	
\usepackage{multirow}
\usepackage{booktabs}
\usepackage{graphicx}
\usepackage{calc}
\usepackage{tikz}
	\usetikzlibrary{calc}
	\usetikzlibrary{positioning}
	\usetikzlibrary{intersections}
	\usetikzlibrary{arrows}
	\usetikzlibrary{backgrounds}
\pgfdeclarelayer{background}
\pgfdeclarelayer{foreground}
\pgfsetlayers{background,main,foreground}   

\usepackage[all,cmtip,2cell]{xy}
\UseTwocells
\usepackage{array}
\usepackage[labelfont=bf]{caption}
\usepackage{tabularx}
	\newcolumntype{s}{>{\hsize=.45\hsize}X}
\allowdisplaybreaks
\newcommand{\colimCC}
{\operatornamewithlimits{\underset{(F\downarrow_{\mathcal{D}}\,d)}{colim}}}
\newcommand{\colimCD}
{\operatornamewithlimits{\underset{(F\downarrow_{\mathcal{D}}\,d^{\prime})}{colim}}}

\usepackage{environ}
\makeatletter
\newsavebox{\measure@tikzpicture}
\NewEnviron{scaletikzpicturetowidth}[1]{%
  \def\tikz@width{#1}%
  \begin{lrbox}{\measure@tikzpicture}%
  \BODY
  \end{lrbox}%
  \pgfmathparse{#1/\wd\measure@tikzpicture}%
  \BODY
}
\makeatother

\newtheorem{definition}{Definition}
\theoremstyle{definition}

\newtheorem{axiom}{Axiom}
\theoremstyle{definition}

\newtheoremstyle{remark}{3pt}{3pt}{}{}{\itshape\bfseries}{}{.5em}{}
\theoremstyle{remark}
\newtheorem{remark}{Remark}[section]

\newtheoremstyle{example}{3pt}{3pt}{}{}{\itshape\bfseries}{}{.5em}{}
\theoremstyle{example}
\newtheorem{example}{Example}[section]

\newcommand{\xdownarrow}[1]{
  {\left\Downarrow\vbox to #1{}\right.\kern-\nulldelimiterspace}
}

\newcommand{\pullback}[2]{{}_{#1}\kern-\scriptspace{\times}_{#2}}

\newcommand{\one}{\textit{One}$^{3}$}
\newcommand{\fourmins}{\textit{4}$^{\prime}$\textit{33}$^{\prime\prime}$}
\newcommand{\zeromins}{\textit{0}$^{\prime}$\textit{00}$^{\prime\prime}$}
\newcommand{\zerominsfull}{\textit{0}$^{\prime}$\textit{00}$^{\prime\prime}\textit{(\fourmins\space No. 2)}$}
\newcommand{\onefull}{\one$=$\fourmins\space (\zeromins) + \raisebox{-1mm}{\includegraphics[width=5pt]{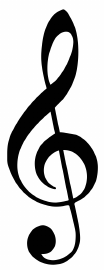}}}

\title{A category-theoretic approach to modeling John Cage's \textit{Silent piece}}
\author{Michael D. Fowler  \\
	Technische Universit\"{a}t Berlin  \\
	}

\date{\today}
\begin{document}

\maketitle
  

\begin{abstract}
\noindent In this article we derive a schema of John Cage's meta-work the \textit{Silent piece} from \fourmins, \zerominsfull, and \onefull, using the mathematics of category theory within Spivak and Kent's (\citeyear{Spivak2012}) framework of ontological logs for knowledge representation. A category presentation $\mathcal{A}$ of a database that describes an instance of \fourmins\space from its premiere in 1952 is translated via two functors into the category presentations $\mathcal{B}$ (\zeromins) and $\mathcal{C}$ (\onefull). A pushout of $\mathcal{B}$ and $\mathcal{C}$ along $\mathcal{A}$ allows for the presentation of the category $\mathcal{S}$ (the meta-work the \textit{Silent piece}), and a discussion of the category's $\mathcal{S}$-specification and fiber order. Finally, we derive a semantics from the fiber in order to reason on persistent spatio-temporal structures of Cage's \textit{Silent piece}.
\end{abstract}

\section{Introduction}
John Cage's \fourmins\space from 1952, \zerominsfull\space from 1962, and \onefull\space from 1989, are considered by the composer as a family of compositions. \citep{Cage2016} Though on the surface, each of the works differ regarding their performance practices and methods for generating sounds, they are united in their approach to problemitizing those fundamental canonical notions of what defines music. The trilogy spans Cage's entire oeuvre, from his pre-\textit{I-Ching} period through to his time bracket approach in the late number pieces.

That this trilogy is strongly connected is most obviously suggested in the titles of each piece, which include successive references to previous compositions. In terms of \one, the last work of the trilogy, Cage describes its origin as emerging from his reluctance to perform \fourmins\space another time:

\begin{quote}
They asked me to perform the silent piece --- a Japanese group was giving a concert in Nagoya --- and it was a concert at the time I was given the Kyoto Prize in Kyoto. I said, I don't want to do the silent piece, because I thought that silence had changed from what it  was, and I wanted to indicate that. \citep[page 94]{Fetterman1996} 
\end{quote}

If we consider that the final work \one\space is Cage's mature conceptualization of his aesthetics of \textit{silence}, then it is clear that there was not only a subtle conceptual shift between \fourmins\space and its second version, \zeromins, but between \zeromins\space and \one. 
In \fourmins, the well-documented premiere \citep{Gann2010} consisted of pianist David Tudor performing a \textit{tacit} over three movements which were indicated by the pianist opening and closing the lid of the piano. Though the work ran for four minutes and thirty-three seconds, Cage later indicates in the 1967 \textit{Second Tacet Edition} score that the work can be for any number of players and any length of time. 

The controversy that the work generates \citep{Dodd2018} emerges from Cage's demand for a new mode of listening and a new musical semiotics.
Brooks' suggests that \fourmins\space can be understood as a \textit{pragmatic test} that emerged after the composer's 1951 visit to an anechoic chamber\footnote{An anecdote that was widely communicated by Cage in his lifetime was the experience of two distinct sounds after entering the Kroft Laboratory's anechoic chamber at Harvard. Cage was told that ``the high one was your nervous system in operation. The low one was your blood circulation.'' See \citet[page 134]{Cage1968}, and \citet[page 64]{Foy2010}.}:
\begin{quote}
If the `silence' of the anechoic chamber was filled with sound, it must have been the case that all previous silences were also filled with sounds---sounds to which Cage had not attended. And since music, too, is filled with sounds, the pragmatic test tells us that any difference between the two (music and silence) cannot depend upon the presence of sound. Sound is not suitable to `verify,' to `validate,' the ideas of music or silence; if sound is to be the criterion, music and silence are identical. \citep[page 125]{Brooks2007} 
\end{quote}

\fourmins\space therefore presents a semiotics that is a de-centering of what \citet{Greimas1983} called the \textit{axis of knowledge} and its intersection with the \textit{axis of desire} within a narrative. The axis of knowledge consists of the triple \{Sender, Object, Receiver\}, and corresponds to the roles of composer, the performative act, and the listener. The axis of knowledge consists of the pair \{Subject, Object\}, and corresponds to the performer and the performative act. \fourmins\space problemitizes the Subject-Receiver pair by proposing a mapping of either the Subject or the Receiver to its neighboring axis. 
This results in the potential of listeners becoming acoustic activators/performers, which is a competing model to standard 19th Century expectations of concert hall etiquette. \citep{Whittington2013,Fowler2019} 

In the follow-up to \fourmins, 1962's \zeromins\space is an about face in which active gestures are employed, though ones that are de-contextualized. Here, the sounds of the piece are a result of the performance of an obligation with a set of restrictions, and using an amplification system set to a maximum level without feedback. The notion of \textit{silence} and the observation of a timeframe from \fourmins\space appears to be mapped to the notion of some non-rareified action and its duration, and what \citet{Kahn1997}, \citet{Pritchett1996}, and \citet{Craenen2014} identify as the extension of the audibility of everyday sounds through theater. In \zeromins, the restrictions include:

\begin{enumerate}
\item{No exact duplication of performances by the person performing the act.}
\item{The act can be subject to interruption(s).}
\item{The act must be performed under conditions of maximum amplification.}
\item{The act cannot be the performance of a work of music.}
\item{There should be no attentiveness paid to the context of the performance: either through electronic, musical or theatrical means.}
\item{The act should be an obligation to others rather than to the person performing the act.}
\end{enumerate}

Though these restrictions are found directly in the score, both \citet[page 120]{Pisaro2018} and \citet[page 196]{Gann2010} suggest the importance of implementing an everyday action in \zeromins, while de-emphasizing the central theme of the piece as a performance of a restricted obligation under maximum amplification.

The use of amplification is also found in the final work in the trilogy, the unpublished work \one. Here, a performance space is maximally amplified, for which the audience and performer simply observe a time frame through intense listening. Cage provides the following account of the premiere of the work in Nagoya in 1989:

\begin{quote}
So what I did was to come on to the stage in front of the audience, and then the feedback level of the auditorium space was brought up to feedback level through the sound-system. There was no actual feed-back, but you knew that you were on the edge of feedback --- which is what I think our environmental situation is now. \citep[page 94]{Fetterman1996} 
\end{quote}

Thus, like \fourmins, its sounds are not connected to an activating musical gesture on an instrument, and like \zeromins, these sounds are made audible through the use of an amplification system. For \citet[page 133]{Swed1993}, this suggests that ``Cage seems to be demonstrating the vague discomfort that unnoticed amplification can have on our sensibilities,'' while \citet{Pritchett2018} identifies the work as ``an example of Cage's darker side,'' and the composer's desire to equate the dangers of highly amplified acoustic spaces to wider environmental and social tipping points that emerged during the 1980s.

In order to ground our investigation into the three works, we use Brooks' notion of a \textit{pragmatics of silence}. Here, the \textit{Silent piece} is framed as Cage's pragmatic test of his earlier anechoic chamber experience through a musical proposal to ``attend to attention'' \citep[page 126]{Brooks2007} via a disciplined action that is viewed as a utility to unlock a world ``teeming with possibilities for experience.'' \citep[page 148]{Johnson1970} For Cage, silence was simply ``the multiplicity of activity that constantly surrounds us. We call it `silence' because it is free of \textit{our} activity.'' \citep[page 245]{KostelanetzCage2003}

\begin{axiom}[Pragmatics of silence]\label{Axiom1} \normalfont
John Cage's musical philosophy of silence is discretely manifested in the works: \fourmins\space(1952), \zeromins\space(1962), and \onefull\space(1989). These works are instances of a meta-work called the \textit{Silent piece} \citep{Cage2016}. For each of the instances we assert the following:
\begin{enumerate}
\item{A sound source, whether intentional or non-intentional, natural or artificial, perceptible or imperceptible, is an element of the set of all possible musical sound sources, and thus satisfies inclusion as a thematic material of an instance of the \textit{Silent piece}.}\label{A1}
\item{A gesture used by a performer in the realization of any instance of the \textit{Silent piece} satisfies the following criteria: 
\begin{enumerate}
\item{The gesture adheres to a given set of local instructions and/or restrictions accompanying an instance.}
\item{The gesture is the performance of an action in a disciplined manner.}
\end{enumerate}
Furthermore, for any instance of the \textit{Silent piece}, there also exists a set of all valid performative actions (gestures) that can be used to realize an instance such that this set is non-empty.}\label{A2}
\end{enumerate}
\end{axiom}

The motivation for this article is to formalize the \textit{pragmatics of silence} through schemas that describe the relations between canonic performance practices, spatio-temporal qualities (sounds and sites of performance), and musical instructions within the three works. 
In order to achieve this end, we apply the mathematics of category theory \citep{Maclane2013}, and Spivak's (\citeyear{Spivak2012}) framework of ontological logs. This approach allows us to derive the meta-work of the \textit{Silent piece} as a categorical schema replete with a fiber order, and an accompanying formal language that describes instances. We use this formal language in order to reason on a number of persistent spatio-temporal structures in the \textit{Silent piece}.
\section{Ologs, databases, and category presentations}

Though category theory has been adept at providing insights into many branches of mathematics as a \textit{meta-language} for understanding algebras of functions \citep{Awodey2010,Maclane2013}, its application into the field of knowledge representation by Spivak has provided a method in which to apply it to real-world contexts. The approach of \citet{Spivak2013b}, and \citet{Spivak2012} focusses on aligning pure mathematics to relational database theory by firstly introducing the notion of a \textit{categorical schema} that consists of a \textit{directed multi-graph} and a set of \textit{path equivalences}, in addition to the labelling of a schema as an \textit{ontological log}. 

\begin{definition}[Directed multi-graph]\label{Graph}\normalfont
A \textit{directed multi-graph} $G$, denoted $G=\{A,V,src,tar\}$, contains of a set of vertices, $V$, and a set of arrows, $A$, for which a source function is denoted $src: A \rightarrow V$, and a target function is denoted $tar: A \rightarrow V$. Given $a \in A$ is an arrow with source $src(a) =v$, and a target $tar(a) =w$, the graph is drawn as:
\begin{equation}
v \xrightarrow{a} w.
\end{equation}
\end{definition}

\begin{definition}[Path length]\label{Pathlength}\normalfont
Given a directed multi-graph $G=\{A,V,src,tar\}$, a \textit{path length n} in $G$ is denoted $p \in \textrm{Path}_{G}^{n}$, and is a head-tail sequence of arrows in $G$:
\begin{equation}
p=(v_{0}\xrightarrow{a_{1}} v_{1}\xrightarrow{a_{2}} \dots \xrightarrow{a_{n}}v_{n}).
\end{equation}
Regarding $\textrm{Path}_{G}^{1}=A$, and $\textrm{Path}_{G}^{0}=V$, we call length 0 on vertex $v$ the \textit{trivial path} on $v$, and denote it $\textrm{id}_{v}$. The set of all paths in $G$ is given through:
\begin{equation}
\textrm{Path}_{G}:=\bigcup_{n \in N} \textrm{Path}_{G}^{n}.
\end{equation}
\end{definition}

\begin{remark}\label{PathRemark}
Every path in a directed multi-graph $G$ contains a source vertex and target vertex, for which given a path $p$ with $src(p) = a$, and $tar(p) = b$, we denote the path as $p: a \rightarrow b$. Furthermore, the set of all paths between $a$ and $b$ is written $\textrm{Path}_{G}(a,b)$. We can also consider a composition of paths given $p: a\rightarrow b$, and $q: b\rightarrow c$, which we denote $q\circ p: a \rightarrow b$.
\end{remark}

\begin{definition}[Categorical path equivalence]\label{PathEq}\normalfont
Given a directed multi-graph $G=\{A,V,src,tar\}$, we define a \textit{categorical path equivalence relation} on $G$ as an equivalence relation $\simeq$ on $\textrm{Path}_{G}$, with the following properties:
\begin{enumerate}\setlength\itemsep{0.1em}
\item If $p \simeq q$, then $src(p) \simeq src(q)$.
\item If $p \simeq q$, then $tar(p) \simeq tar(q)$.
\item Given two paths $p,q: b \rightarrow c$, and the arrow $m: a\rightarrow b$, if $p \simeq q$, then $mp \simeq mq$.
\item Given the two paths $p,q: a \rightarrow b$, and the arrow $n: b\rightarrow c$, if $q\simeq q$, then $pn \simeq qn$.
\end{enumerate}
\end{definition}

\begin{definition}[Categorical schema]\label{CatSch}\normalfont
A \textit{categorical schema} is a 2-tuple $\mathcal{C}=(G, \simeq)$, in which $G$ is a directed multi-graph, and $\simeq$ is a categorical path equivalence on $G$. An \textit{ontological log} (or simply \textit{olog}) of a categorical schema is the labelling of the vertices $V$, and arrows $A$ of $G$ such that:
\begin{itemize}\setlength\itemsep{0.1em}
\item a vertex $v$ represents a \textit{type} or set of objects that form a class, for which its label is written as a \textit{singular indefinite noun phrase} inscribed in a rectangle.
\item an arrow $a$ represents an \textit{aspect} or function between types, for which its label is written as a verb phrase.
\item an olog is accompanied with a set of equations that are path equivalences statements
\end{itemize}
\end{definition}

\begin{remark}
It is perhaps self-evident that Spviak's definition of a categorical schema is equivalent to a standard definition of a category \citep{Maclane2013}. For example, the set of vertices $V_{\mathcal{C}}$ in a categorical schema $\mathcal{C}$ are the \textit{objects} found in $\textrm{Obj}(\mathcal{C})$, while the arrows $A_{\mathcal{C}}$ are the \textit{morphisms} found in $\textrm{Hom}_{\mathcal{C}}$, that have objects that are sources and targets. Similarly, the \textit{identity morphism} $\textrm{id}_{X}$ on object $X$ is the trivial path, and the composition and associativity laws are covered through a schema's congruence rule and path equivalence. Furthermore a categorical schema is a small category given that type labels indicate sets of things. In this article we use Spivak's terminology interchangeably with standard terms in category theory.
\end{remark}

\begin{remark}
The choice of a type is a distinction made by the author of an olog, for which its label describes its \textit{ontological intention}, and for which its set of instances represent its \textit{ontological extension}. Instances of a type can thus be useful in documenting a type's functional relation to other types. Therefore, a label of an aspect on an arrow $a$ refers to a \textit{functional dependence} between a source type, or domain, and target type, or co-domain. That is, each instance of a source type should have exactly one corresponding target type. An olog then enables the parsing of sentences in natural language. For example the following directed multi-graph: 
\begin{equation}
\begin{tikzpicture}[thick]\small
\tikzstyle{D} = [draw,fill=none,align=center,inner sep=6]
\node[D,fill=white] at (0,0) (A) {a musical score};
\node[D] at (6,0) (C) {a person};
\draw[->] (A) --node[above,align=center]{is composed by}(C);
\end{tikzpicture}
\end{equation}
can be parsed as the sentence: ``A musical score is composed by a person.'' Given the functional nature of the aspect ``is composed of,'' the graph asserts a \textit{world view} that: 
\begin{enumerate}
\item{There may be many musical scores composed by a single person.} 
\item{There are no musical scores that do not have a person as author.}
\item{There are no musical scores composed by more than one person.}
\end{enumerate} 
The functional nature of aspects between types in an olog also means that more expressive types can be constructed from sequences such as the \textit{cartesian product} of types.
\end{remark}

\begin{definition}[Product]\label{Product}\normalfont
Given a categorical schema $\mathcal{C}$, whose vertices (types) contain the sets $X$ and $Y$, we denote their \emph{product} as $X\times Y$, such that the set of ordered pairs $(x,y)$ can be expressed as:
\begin{equation}
X \times Y = \{(x,y)\,|\,x \in X \wedge y \in Y\}.
\end{equation}
There also exists two functions called \textit{projection functions}, $\tau_{1}$ and $\tau_{2}$, that share a domain:
\begin{equation}
\begin{tikzpicture}
\node at (0,0)(T) {$X \times Y$};
\node at (-1.5,-1.5)(X){$X$};
\node at (1.5,-1.5)(Y){$Y$};
\draw[->](T)--node[left]{$\tau_{1}$}(X);
\draw[->](T)--node[right]{$\tau_{2}$}(Y);
\end{tikzpicture}
\end{equation}
We label the olog of $\mathcal{C}$ containing the product $X\times Y$ of two vertices as:
\begin{equation}
``X\times Y":=\textrm{a pair }(a,b) \textrm{ where } a \textrm{ is } ``X" \textrm{ and } b \textrm{ is } ``Y,"
\end{equation}
for which the projection functions $\tau_{1}$, and $\tau_{2}$ are labelled as ``yields as $a$,'' and ``yields as $b$'' respectively. \citet[page 37]{Spivak2013} There also exists a \textit{universal property of products} in which given a vertex (set) $A$, the functions $f_{1}:A\rightarrow X$ and $f_{2}:A\rightarrow Y$ induce a unique function, $f: A\rightarrow X\times Y$, such that $\tau_{1} \circ f \simeq f_{1}$, and $\tau_{2} \circ f \simeq f_{2}$:
\begin{equation}
\begin{tikzpicture}
\node at (0,0)(T) {$X \times Y$};
\node at (-1.5,-1.5)(X){$X$};
\node at (1.5,-1.5)(Y){$Y$};
\node at (0,-3)(A){$A$};
\draw[->](T)--node[left]{$\tau_{1}$}(X);
\draw[->](T)--node[right]{$\tau_{2}$}(Y);
\draw[->,dashed](A)--node[left,pos=0.5]{$\exists !f$}(T);
\draw[->](A)--node[left,pos=0.5]{$f_{1}$}(X);
\draw[->](A)--node[right,pos=0.5]{$f_{2}$}(Y);
\end{tikzpicture}
\end{equation}
We follow the suggestion of \citet[page 37]{Spivak2013}, and label the unique arrow $f$ in an olog as:
\begin{equation}
``f": \textrm{ yields, insofar as it } ``f_{1}"\,\, ``X" \textrm{ and } ``f_{2}"\,\, ``Y."
\end{equation}
\end{definition}

But in order to utilize the migration of information contained within categorical schemas, Spivak introduces the framework of the \textit{translation} of schemas. \citep[page 7]{Spivak2013b}

\begin{definition}[Translation of schemas]\label{Trans}\normalfont
Let $G=\{A_{G},V_{G},src_{G},tar_{G}\}$, and $H=\{A_{H},V_{H},src_{H},tar_{H}\}$ be two directed multi-graphs such that $\mathcal{C}=(G,\simeq_{\mathcal{C}})$, and $\mathcal{D}=(H,\simeq_{\mathcal{D}})$ are their associated categorical schemas. A \textit{translation}, denoted $F:\mathcal{C}\rightarrow\mathcal{D}$ contains the following elements:
\begin{enumerate}\setlength\itemsep{0.1em}
\item a function $V_{F}: V_{G}\rightarrow V_{H}$,
\item a function $A_{F}:A_{G}\rightarrow \textrm{Path}_{H}$,
\end{enumerate}
Additionally, the function $A_{F}$ preserves both paths equivalences as well as targets and sources such that the following diagrams commute:
\begin{equation}
\begin{tikzpicture}[scale=2,thick]
\node at (0,0) (AG) {$A_{G}$};
\node at (1,0) (PathH) {$\textrm{Path}_{H}$};
\node at (0,-1) (VG) {$V_{G}$};
\node at (1,-1) (VH) {$V_{H}$};
\draw [->](AG)--node[above,pos=0.5]{$A_{F}$}(PathH);
\draw [->](AG)--node[left,pos=0.5]{$src_{G}$}(VG);
\draw [->](PathH)--node[right,pos=0.5]{$src_{H}$}(VH);
\draw [->](VG)--node[below,pos=0.5]{$V_{F}$}(VH);
\end{tikzpicture}\hspace{5mm}
\begin{tikzpicture}[scale=2,thick]
\node at (0,0) (AG) {$A_{G}$};
\node at (1,0) (PathH) {$\textrm{Path}_{H}$};
\node at (0,-1) (VG) {$V_{G}$};
\node at (1,-1) (VH) {$V_{H}$};
\draw [->](AG)--node[above,pos=0.5]{$A_{F}$}(PathH);
\draw [->](AG)--node[left,pos=0.5]{$tar_{G}$}(VG);
\draw [->](PathH)--node[right,pos=0.5]{$tar_{H}$}(VH);
\draw [->](VG)--node[below,pos=0.5]{$V_{F}$}(VH);
\end{tikzpicture}
\end{equation}
\end{definition}

We see that translations of schemas are actually \textit{functors}, given that they send vertices to vertices, arrows to arrows, and path equivalences to path equivalences. They can thus also be considered graph morphisms on $G$. By this we mean that given two directed multi-graphs $G$ and $G^{\prime}$, a functor $F$ is the \textit{graph homomorphism} $F:G\rightarrow \textrm{Path}_{G^{\prime}}$, in which for any paths $p$ and $q$, if $p\simeq q$ in $G$, then $\textrm{Path}_{F}(p)\simeq_{G^{\prime}} \textrm{Path}_{F}(q)$.

This allows for the definition of the \textit{category of categorical schemas}, $\mathbf{Sch}$, in which objects are categorical schemas (categories), and morphisms are translations (functors). \citet{Spivak2013b} proves $\mathbf{Sch}\simeq \mathbf{Cat}$ through constructing a category via functors from a schema and schema from a category, for which we can understand $\mathbf{Sch}$ as differentiated from $\mathbf{Cat}$ through the idea that $\mathcal{C}$ is a \textit{category presentation}. It is a presentation given that it is governed by objects and arrows that are generators, and path congruences that are like relations.

In terms of the generation of \textit{instances} on a schema, the functor given by $I: \mathcal{C}\rightarrow \mathbf{Set}$, where $\mathbf{Set}$ is the category of sets whose objects are sets and morphisms are functions, allows us to consider a schema as a representation of a database in which objects are tables (with IDs) that contain values, and arrows are columns pointing towards other tables. 

\begin{definition}[Instances on $\mathcal{C}$]\label{Inst}\normalfont
Given a categorical schema, $\mathcal{C}:=(G,\simeq)$, such that $G=(A,V,src,tar)$, we call the functor $I$ an \textit{instance} on $\mathcal{C}$ such that $I$ has following properties:
\begin{enumerate}\setlength\itemsep{0.1em}
\item $\forall v \in V$, there exists a set $I(v)$.
\item $\forall a \in A$, where $a: v\rightarrow v^{\prime}$, the exists a function $I(a): I(v) \rightarrow I(v^{\prime})$.
\item For every path equivalence $p\simeq q$, $I(p) = I(q)$ is true.
\end{enumerate}
\end{definition}

\begin{example}
Consider the following olog (categorical schema) that has a path equivalence which states that there exists musical works that are parts of other musical works that are composed by some person.

\begin{align}
\begin{tikzpicture}[thick]\small
\tikzstyle{D} = [draw,fill=none,align=center,inner sep=6]
\node[D,fill=white] at (0,0) (A) {a musical work};
\node[D] at (6,0) (C) {a person};
\draw[->] (A) --node[above,align=center]{is composed by}(C);
\path[->] (A) edge [loop above] node {is part of} ();
\end{tikzpicture}\\[2.2ex]
(\textrm{is composed by})\circ(\textrm{is part of})\simeq (\textrm{is composed by})
\end{align}
As a database of instances in which information can be collated under labelled type vertices, we have the following tables given through the functor $I: \mathcal{C}\rightarrow \mathbf{Set}$:
\small{
\begin{align}
\begin{tabularx}{0.9\textwidth}{| l || >{\raggedright\arraybackslash}X | l | >{\raggedright\arraybackslash}X |}
\hline
\multicolumn{3}{|c|}{\cellcolor{white}\textbf{a musical work}}\\\hline
\textbf{ID}&\textbf{is part of}&\textbf{is composed by}\\\hline
\textit{December 1952}&\textit{Folio}&Earle Brown\\\hline
\textit{Oktophonie}&\textit{Dienstag aus Licht}&Karlheinz Stockhausen\\\hline
\end{tabularx}
\nonumber\\[1.2ex]
\begin{tabularx}{0.4\textwidth}{| >{\raggedright\arraybackslash}X ||}
\hline
\multicolumn{1}{|c|}{\cellcolor{white}\textbf{a person}}\\\hline
\textbf{ID}\\\hline
Earle Brown\\\hline
Karlheinz Stockhausen\\\hline
\end{tabularx}
\end{align}
}
\end{example}

IDs track instances along rows of an object in the categorical schema, and columns point to other tables, and thus track functions between objects. Columns are also known as \textit{foreign keys} given they cross-reference data.

\section{\fourmins\space as categorical schema $\mathcal{A}$}
 
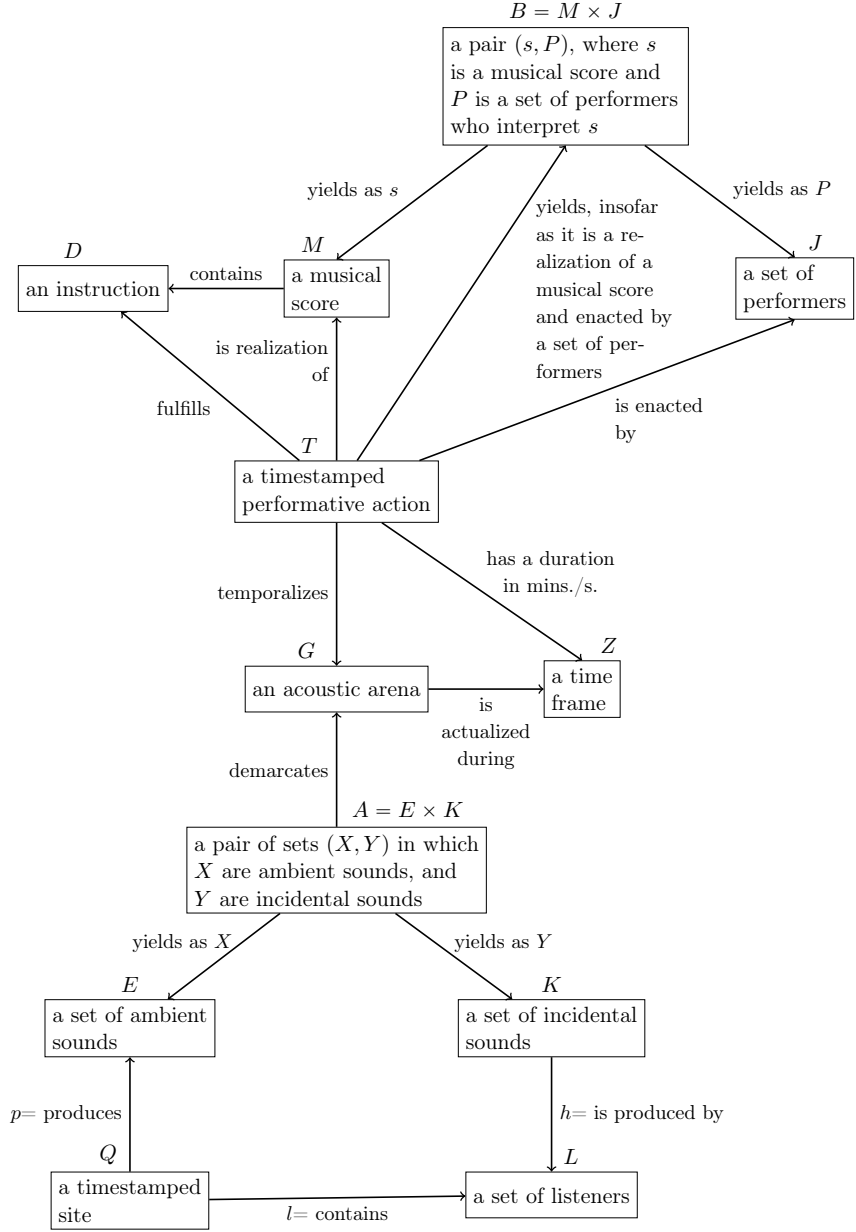
\begin{figure}[th!]
\centering
\adjustbox{max width=0.9\textwidth}{
\begin{tikzpicture}
	\tikzset{
		V-node/.style={
		fill=white,draw=black,draw,minimum 	size=8mm
		},
}
\node[V-node,label={[align=left]97:$M$},align=left] at (0,0) (MusicalScore) {a musical\\score};
\node[V-node,right=60mm of MusicalScore,align=left,label={[align=left]80:$J$},align=left](Sequence){a set of\\ performers};
\node[V-node,left = 20mm of MusicalScore,align=left,label={[align=left]97:$D$}](Instruction){an instruction};
\node[V-node,below=25mm of MusicalScore,align=left,label={[align=left]110:$T$}] (PerformativeAction){a timestamped\\performative action}; 
\node[V-node,below=25mm of PerformativeAction,align=left,label={[align=left]120:$G$}](AcousticArena){an acoustic arena}; 
\node[V-node,xshift=0,label={[align=left]70:$Z$},align=left,right=20mm of AcousticArena] (TimeFrame){a time\\ frame};
\node[V-node,align=left,below=20mm of AcousticArena,label={[align=left]80:$A=E\times K$}](Pair){a pair of sets $(X, Y)$ in which\\$X$ are ambient sounds, and\\ $Y$ are incidental sounds};
\node[V-node,align=left,below left=15mm and -5mm of Pair,label={[align=left]90:$E$}](AmbientSounds){a set of ambient\\sounds};
\node[V-node,align=left,below right=15mm and -5mm of Pair,label={[align=left]90:$K$}] (IncidentalSounds){a set of incidental\\sounds};
\node[V-node,above=25mm of PerformativeAction,label={[align=left]90:$B= M\times J$},align=left,text width=40mm] at ($(MusicalScore)!0.5!(Sequence)$) (Performer){a pair $(s,P)$, where $s$ is a musical score and $P$ is a set of performers who interpret $s$};
\node[V-node,below=20mm of AmbientSounds,label={[align=left]100:$Q$},align=left](Site){a timestamped\\site};
\node[V-node,below=20mm of IncidentalSounds,label={[align=left]80:$L$},align=left](Listener){a set of listeners};
\draw[thick,<-](TimeFrame)--node[below,pos=0.5,align=center]{\small is\\\small actualized\\\small during}(AcousticArena);
\draw[thick,->](PerformativeAction)--node[pos=0.6,above,align=center,xshift=24]{\small has a duration\\\small in mins./s.}(TimeFrame.north);
\draw[thick,->](MusicalScore)--node[above,pos=0.5]{\small contains}(Instruction);
\draw[thick,<-](MusicalScore.south)--node[left,pos=0.3,align=center,align=right,xshift=0]{\small is realization\\\small of}(PerformativeAction);
\draw[thick,->](PerformativeAction)--node[left,pos=0.5,align=center]{\small temporalizes}(AcousticArena);
\draw[thick,->](Pair)--node[left,pos=0.5,align=center]{\small demarcates}(AcousticArena);
\draw[thick,->](Pair)--node[left,pos=0.35,align=center]{\small yields as $X$}(AmbientSounds);
\draw[thick,->](Pair)--node[right,pos=0.35,align=center,xshift=5]{\small yields as $Y$}(IncidentalSounds);
\draw[thick,->](Site)--node[left,pos=0.5,align=center]{\small $p$= produces}(AmbientSounds);
\draw[thick,<-](Listener)--node[right,pos=0.5,align=center]{\small $h$= is produced by}(IncidentalSounds);
\draw[thick,<-](Listener)--node[below,pos=0.5,align=center]{\small $l$= contains}(Site);
\draw[thick,->](Performer)--node[right,pos=0.4,align=center,xshift=10]{\small yields as $P$}(Sequence.north);
\draw[thick,->](Performer)--node[left,pos=0.4,align=center,xshift=-10]{\small yields as $s$}(MusicalScore.north);
\draw[thick,<-](Sequence.south)--node[pos=0.5,yshift=-15,right,align=left]{\small is enacted\\\small by}(PerformativeAction);
\draw[thick,->](PerformativeAction)--node[left,pos=0.35,align=center,xshift=-9]{\small fulfills}(Instruction);
\draw[thick,->](PerformativeAction)--node[right,pos=0.45,align=left,yshift=15,text width=25mm,xshift=38]{\small yields, insofar as it is a realization of a musical score and enacted by a set of performers}(Performer.south);
\end{tikzpicture}
}
\caption{Olog of $\mathcal{A}$.}
\label{OlogA}
\end{figure}

Consider the olog of the categorical schema $\mathcal{A}$ found in Figure~\ref{OlogA} as a representation of the background knowledge about the premiere of \fourmins\space in Woodstock, NY 1952. The olog presents not only relations between agents and objects, but also spatial and acoustic conditions. In particular, the type `an acoustic arena' ($G$) is what \citet{BlesserSalter2007} define as an auditory space in which sound sources can be heard by an acoustic community defined over a physical coverage area. A concert hall is a typical manifestation of an acoustic arena discretely defined across an architectural interior, though a town square or stadium are also associated to acoustic arenas. An acoustic arena can be interior to some other arena or arenas, and therefore nested, but acoustically disjoint from its neighbors. 

The olog contains the type (object) $A=E\times K$ that is a pair of sets of ambient sounds and incidental sounds, of which the ambient sounds are produced at, or by the site of the performance. Such sounds were confirmed by Tudor \citep[page 4]{Gann2010} during the premiere, for which incidental sounds came from confused audience members during the performance. We also declare the olog $\mathcal{A}$ as containing the following path equivalences:
 \begin{align}
&(\textrm{fulfills})\simeq(\textrm{contains})\circ(\textrm{is realization of})\label{eqv1}\\[1em]
&(\textrm{has duration in mins/s})\simeq(\textrm{is actualized during})\circ(\textrm{temporalizes})\label{eqv2}\\[1em]
&(\textrm{yields as } s) \circ (\textrm{yields, insofar as it is a realization of a musical score }\nonumber\\
&\textrm{and enacted by a set performers})\simeq(\textrm{is realization of})\label{eqv3}\\[1em]
&(\textrm{yields as } P) \circ (\textrm{yields, insofar as it is a realization of a musical score }\nonumber\\
&\textrm{and enacted by a set performers})\simeq(\textrm{is enacted by})\label{eqv4}
 \end{align} 
These equivalences assert that in equation~\eqref{eqv1}, a timestamped performative action fulfills an instruction, given that the instruction is part of the musical work being realized. This path equivalence points to an assessment of \textit{authenticity} about a realization of the work. Regarding the second case in equation~\eqref{eqv2}, we can assert that the timestamped performative action has a duration that is equal to the amount of time by which some acoustic arena is actualized. This elides with what \citet{Davies2003}, \citet{Shultis2013}, and \citet{Whittington2013} acknowledge as a process of \textit{temporal framing} underway in the work. Here, the notion of a passive gesture of not playing an instrument, or what Retallack considers a ```non-performance' art [that] is fundamentally performative,'' \cite[page 245]{Retallack2013} 
results in a mapping of the pair $A=E\times K$ of incidental and ambient sounds to an acoustic arena that satisfies point (1) in Axiom~\ref{A1}.

The third and fourth equivalences in equation~\ref{eqv3} and equation~\ref{eqv4} arise from the universal property of products and the object $B=M\times J$. Here the olog points towards the \textit{agency} of the set of performers, for which the timestamped performative action can be considered as evoking an instance of a \textit{performance practice}. 

We note that although the premiere of \fourmins\space was given by pianist Tudor, Cage's later published scores noted that the work could be realized by any number of performers. \citep{Cage1986} The categorical schema $\mathcal{A}$ is based on those conditions of the premiere of \fourmins, but by using the pair $B=(s,P)$, where $s$ is a musical score (\fourmins), and $P$ are a set of performers, we can both account for Tudor's premiere of the work as a solo in 1952 (that is, $P$ as a singleton set), with other accounts of later performances of the piece using various sized ensembles. \citep[page 80-81]{Fetterman1996}
 
\begin{figure}
\scriptsize{
\begin{tabularx}{0.5\textwidth}{| l || >{\raggedright\arraybackslash}X | l |}
\hline
\multicolumn{2}{|c|}{\cellcolor{white}\textbf{a musical score}}\\\hline
\textbf{ID}&\textbf{contains}\\\hline
\fourmins&For any number of performers in three parts (I, II, III) in which each part consists of the performance of a tacit of an agreed length of time\\\hline
\end{tabularx}
}
\scriptsize{
\begin{tabularx}{0.48\textwidth}{| >{\raggedright\arraybackslash}X ||}
\hline
\multicolumn{1}{|c|}{\cellcolor{white}\textbf{an instruction}}\\\hline
\textbf{ID}\\\hline
For any number of performers in three parts (I, II, III) in which each part consists of the performance of a tacit of an agreed length of time\\\hline
\end{tabularx}
}\smallskip

\scriptsize{
\begin{tabularx}{0.7\textwidth}{| >{\raggedright\arraybackslash}X || p{2cm} | p{2cm} |}
\hline
\multicolumn{3}{|p{7cm}|}{\textbf{a pair $(s,P)$, where $s$ is a musical score and $P$ is a set of performers who interpret $s$}}\\\hline
\textbf{ID}&\textbf{yields as $s$}&\textbf{yields as $P$}\\\hline
(\fourmins,\{David Tudor\})&\fourmins&\{David Tudor\}\\\hline
\end{tabularx}
}
\scriptsize{
\begin{tabularx}{0.25\textwidth}{|>{\raggedright\arraybackslash}X ||}
\hline
\multicolumn{1}{|>{\raggedright\arraybackslash}X|}{\cellcolor{white}\textbf{a set of performers}}\\\hline
\textbf{ID}\\\hline
\{David Tudor\}\\\hline
\end{tabularx}
}\smallskip

\scriptsize{
\begin{tabularx}{\linewidth}{
	|>{\raggedright\arraybackslash\hsize=0.35\hsize}X ||
	>{\raggedright\arraybackslash\hsize=0.5\hsize}X |
	>{\raggedright\arraybackslash\hsize=0.3\hsize}X |
	>{\raggedright\arraybackslash\hsize=0.5\hsize}X |
	>{\raggedright\arraybackslash\hsize=0.4\hsize}X |
	>{\raggedright\arraybackslash\hsize=0.29\hsize}X |
	>{\raggedright\arraybackslash\hsize=0.3\hsize}X |
}                              
\hline
\multicolumn{7}{|c|}{\cellcolor{white}\textbf{a timestamped performative action}}\\\hline
\textbf{ID}&\textbf{fulfills}&\textbf{is realization of}&\textbf{yields, insofar as it is a realization of a musical score and enacted by a set of performers}&\parbox[t]{2cm}{\textbf{tempora-}\\\textbf{lizes}}&\textbf{has duration in mins./s.}&\textbf{is enacted by}\\\hline
The observation of a time frame in silence and without intentional musical actions on 29.08.1952&For any number of performers in three parts (I, II, III) in which each part consists of the performance of a tacit of an agreed length of time&\fourmins&(\fourmins, \{David Tudor\})&The interior and immediate surrounds of Maverick Concert Hall&4 minutes and 33 seconds&\{David Tudor\}\\\hline
\end{tabularx}
}\smallskip

\scriptsize{
\begin{tabularx}{0.67\textwidth}{| p{4.8cm} || p{2.7cm}  |}
\hline
\multicolumn{2}{|c|}{\cellcolor{white}\textbf{an acoustic arena}}\\\hline
\textbf{ID}&\textbf{is actualized during}\\\hline
The interior and immediate surrounds of Maverick Concert Hall& 4 minutes and 33 seconds\\\hline
\end{tabularx}
}
\scriptsize{
\begin{tabularx}{0.3\textwidth}{| >{\raggedright\arraybackslash}X ||}
\hline
\multicolumn{1}{|c|}{\cellcolor{white}\textbf{a timeframe}}\\\hline
\textbf{ID}\\\hline
4 minutes and 33 seconds\\\hline
\end{tabularx}
}\smallskip

\scriptsize{
\begin{tabularx}{\textwidth}{| >{\hsize=0.25\hsize\raggedright\arraybackslash}X || >{\hsize=0.25\hsize\raggedright\arraybackslash}X | >{\hsize=0.25\hsize\raggedright\arraybackslash}X | >{\hsize=0.2\hsize\raggedright\arraybackslash}X |}
\hline
\multicolumn{4}{|c|}{\cellcolor{white}\textbf{A pair of sets $(X, Y)$ in which $X$ are ambient sounds, and $Y$ are incidental sounds}}\\\hline
\textbf{ID}&\textbf{demarcates}&\textbf{yields as $X$}&\textbf{yields as $Y$}\\\hline
(\{wind in trees, raindrops on roof\}, \{talking, rustling paper, walking-out\})&The interior and immediate surrounds of Maverick Concert Hall&\{wind in trees, raindrops on roof\}&\{talking, rustling paper, walking-out\}\\\hline
\end{tabularx}
}\smallskip

\scriptsize{
\begin{tabularx}{0.5\textwidth}{| >{\raggedright\arraybackslash}X ||}
\hline
\multicolumn{1}{|c|}{\cellcolor{white}\textbf{a set of ambient sounds}}\\\hline
\textbf{ID}\\\hline
\{wind in trees, raindrops on roof\}\\\hline
\end{tabularx}
}
\scriptsize{
\begin{tabularx}{0.48\textwidth}{| >{\raggedright\arraybackslash}X || >{\raggedright\arraybackslash}X |}
\hline
\multicolumn{2}{|c|}{\cellcolor{white}\textbf{a set of incidental sounds}}\\\hline
\textbf{ID}&\textbf{is produced by}\\\hline
\{talking, rustling paper, walking-out\}&\{N.N.\}\\\hline
\end{tabularx}
}\smallskip

\scriptsize{
\begin{tabularx}{0.7\textwidth}{| >{\raggedright\arraybackslash}X || >{\raggedright\arraybackslash}X | >{\raggedright\arraybackslash}X |}
\hline
\multicolumn{3}{|c|}{\textbf{A timestamped site}}\\ \hline
\textbf{ID}&\textbf{contains}&\textbf{produces}\\\hline
Maverick Concert Hall, Woodstock, NY, 29.08.1952&\{N.N.\}&\{wind in trees, raindrops on roof\}\\\hline
\end{tabularx}
}
\scriptsize{
\begin{tabularx}{0.238\textwidth}{| >{\hsize=0.25\hsize\raggedright\arraybackslash}X ||}
\hline
\multicolumn{1}{|c|}{\textbf{A set of listeners}}\\ \hline
\textbf{ID}\\\hline
\{N.N.\}\\\hline
\end{tabularx}
}

\caption{The database $\protect\mathbb{D}_{\mathcal{A}}$.}
\label{DB1}
\end{figure}

By applying the functor $I: \mathcal{A}\rightarrow \mathbf{Set}$, we can generate instances, as found in the database $\mathbb{D}_{\mathcal{A}}$ in Figure~\ref{DB1}. The database contains the particulars of the site of the performance, the coverage area of the acoustic arena, and instances of the ambient and incidental sounds that were present at the premiere in 1952. We only include a singleton set of listeners, though as \citet{Gann2010} and \citet{Gentry2018} report, there were possibly a number of audience members who left the concert during the performance in protest.

We could further imagine another functor that provides data about other realizations, such as Cage's 1970 open-air performance of \fourmins\space at Harvard Square. \citep{Katschthaler2016} In this case, we can define a \textit{category of finite instances} on a categorical schema $\mathcal{C}$, for which there exist natural transformations between instances that account for \textit{progressive updates} on a database. 
 
\begin{definition}[Category of finite instances]\label{CatInst}\normalfont
Let $\mathcal{C}$ be a categorical schema, for which the \textit{category of finite instances} on $\mathcal{C}$, denoted $\mathcal{C}$--$\mathbf{Fin}$, contains objects that are instances (functors as $I:\mathcal{C}\rightarrow\mathbf{Set}$), and morphisms that are natural transformations between instances. Given the instances $I,J$ on $\mathcal{C}$, and the morphism $\eta: I \rightarrow J$, we say that for each object $v$, there is a \textit{component} of $\eta$ at $v$, such that the following diagram commutes: 
\begin{equation}
\begin{tikzpicture}[rotate=270,yscale=-1]
\node at (0,0) (FA) {$I(v)$};
\node at (3,0) (FB) {$I(u)$};
\node at (0,-3) (GA) {$J(v^{\prime})$};
\node at (3,-3) (GB) {$J(u^{\prime})$};
\draw[->] (FA)--node[above]{$\eta_{v}$}(GA);
\draw[->] (FB)--node[below]{$\eta_{u}$}(GB);
\draw[->] (FA)--node[left]{$I(a)$}(FB);
\draw[->] (GA)--node[right]{$J(a^{\prime})$}(GB);
\end{tikzpicture}
\end{equation}

We call the natural transformation, $\eta: I\rightarrow J$, a \textit{progressive update} on $I$, which equates to the insertion of rows in a table in $\mathcal{C}$.
\end{definition}

This allows us to also consider a graph of collected instances of performances through a \textit{Grothendieck construction}, and the \textit{discrete opfibration}, $(\int_{\mathcal{C}}(I),\pi_{I})$, where $\int_{\mathcal{C}}(I)$ is the \textit{category of elements} of $I$, and $\pi_{I}$ is a functor that maps instances back onto $\mathcal{C}$. \citep[page 12]{Spivak2013c}

\begin{definition}[Grothendieck construction]\label{G-construct}\normalfont
The \textit{Grothendieck construction} associates the functor $I:\mathcal{C}\rightarrow \mathbf{Set}$ on $\mathcal{C}$ to the pair $(\int_{\mathcal{C}}(I),\pi_{I})$, where $\int_{\mathcal{C}}(I) \in \mathbf{Cat}$, and is the \textit{category of elements} of $I$, and $\pi_{I}$ is the \textit{projection functor} from $\int_{\mathcal{C}}(I)\rightarrow \mathcal{C}$. The construction consists of the following objects and morphisms:
\begin{align}
\textrm{Obj}(\int_{\mathcal{C}}(I)):&=\{(s,x)\,|\,s \in \textrm{Obj}(\mathcal{C}), x \in I(s)\},\\
\textrm{Hom}_{\int_{\mathcal{C}}(I)}((s,x),(s^{\prime},x^{\prime}):&=\{f:s\rightarrow s^{\prime}\,|\,I(f)(x)=x^{\prime}\},
\end{align}
for which the functor $\pi_{I}:\int_{\mathcal{C}}(I)\rightarrow \mathcal{C}$ sends the object $(s,x)$ to $s$ and the morphism $f:(s,x)\rightarrow(s^{\prime},x^{\prime})$ to $f:s\rightarrow s^{\prime}$. 
\end{definition}
 
\begin{figure}[t!]
\centering
\adjustbox{max width=\textwidth}{
\begin{tikzpicture}[every node/.style={font=\scriptsize}]
	\tikzset{
		V-node/.style={
		fill=black,draw=black,draw,minimum 	size=1mm,inner sep=0,shape=circle
		},
}
\node[V-node,label={[align=left]0: \fourmins},align=left] at (0,0) (MusicalScore) {};
\node[V-node,above right=20mm and 15mm of MusicalScore,label={[align=left,yshift=2]180:(\fourmins, \{John Cage\})},align=left] (MusicalScoreString) {};
\node[V-node,right=35mm of MusicalScore,align=left,label={[align=left]0:\{David Tudor\}}](Sequence){};
\node[V-node,right=30mm of MusicalScore,align=left,label={[align=left,xshift=-3,yshift=2]45:\{John Cage\}}](Sequence2){};
\node[V-node,left = 20mm of MusicalScore,align=left,label={[align=left]97:}](Instruction){ };
\node[text width=50mm,inner sep=1,align=left,anchor=south,yshift=3,xshift=-30] at (Instruction){For any number of performers in three parts (I, II, III) in which each part consists of the performance of a tacit of an agreed length of time};
\node[V-node,below left=20mm and 2mm of MusicalScore,align=left,label={[align=left,yshift=4]100:}](PerformativeAction){}; 
\node[V-node,below right=20mm and 2mm of MusicalScore,align=left,label={[align=left]0:}](PerformativeAction2){}; 
\node[text width=50mm,inner sep=1,align=left,anchor=east,xshift=-12,yshift=-8] at (PerformativeAction){The observation of a time frame in silence and without intentional musical actions on ?.?.1970};
\node[text width=50mm,inner sep=1,align=left,anchor=west,xshift=15,yshift=-10] at (PerformativeAction2){The observation of a time frame in silence and without intentional musical actions on 29.08.1952};
\node[V-node,below=20mm of PerformativeAction,align=left,label={[align=left,yshift=-3]120:}](AcousticArena){}; 
\node[V-node,below=20mm of PerformativeAction2,align=left,label={[align=left]45:}](AcousticArena2){}; 
\node[text width=30mm,inner sep=1,align=left,anchor=north,xshift=50,yshift=-13] at (AcousticArena2) {The interior and immediate surrounds of Maverick Concert Hall}; 
\node[text width=30mm,inner sep=1,align=left,anchor=east,xshift=-15,yshift=-5] at (AcousticArena2) {Harvard Square and neighboring streets}; 
\node[V-node,xshift=0,label={[align=left]1:4 minutes and 33 seconds},align=left] at (AcousticArena-|Sequence)(TimeFrame){};
\node[V-node,align=left,below=20mm of AcousticArena,label={[align=left,yshift=-19,text width=40mm]120:(\{wind, traffic, pedestrians\}, \{talking, photography\})}](Pair){};
\node[V-node,align=left,label={[align=left,text width=50mm,yshift=-10,xshift=18]0:(\{wind in trees, raindrops on roof\}, \{talking, rustling paper, walking-out\})}] at (PerformativeAction2|-Pair) (Pair2){};
\node[V-node,align=left,below=20mm of Pair,label={[align=left,text width=20mm,xshift=0,yshift=-20]0:\{wind in trees, raindrops on roof\}},xshift=-10mm](AmbientSounds){};
\node[V-node,align=left,left=4mm of AmbientSounds,label={[align=left,yshift=-10,xshift=-25,text width=20mm]90:\{wind, traffic, pedestrians\}}](AmbientSounds2){};
\node[V-node,align=left,below=20mm of Pair2,label={[align=left,text width=18mm,xshift=0,yshift=-5]160:\{talking, photography\}},xshift=10mm] (IncidentalSounds){};
\node[V-node,align=left,right=4mm of IncidentalSounds,label={[align=left,yshift=-3,xshift=0,text width=27mm]0:\{talking, rustling paper, walking-out\}}] (IncidentalSounds2){};
\node[V-node,right=4mm of MusicalScoreString,label={[align=left,yshift=2]0:(\fourmins, \{David Tudor\})}] (Performer){};
\node[V-node,below=20mm of AmbientSounds,label={[align=left,yshift=-5]-90:}](Site){};
\node[V-node,label={[align=left,yshift=3]100:}] at (AmbientSounds2|-Site)(Site2){};
\node[text width=27mm,inner sep=1,align=left,anchor=north,xshift=35,yshift=-10] at (Site){Maverick Concert Hall, Woodstock, NY, 29.08.1952};
\node[text width=30mm,inner sep=1,align=right,anchor=east,xshift=-5,yshift=-5] at (Site2){Harvard Square, Cambridge, MA, ?.?.1970};
\node[V-node,below=20mm of IncidentalSounds,label={[align=left,yshift=-4,xshift=-5]-80:\{Nam June Paik\}},align=left](Listener){};
\node[V-node,label={[align=left,yshift=2]0:\{N.N.\}},align=left] at (IncidentalSounds2|-Listener)(Listener2){};
\draw[<-](TimeFrame) to [bend right=12]node[below,pos=0.5,align=center]{}(AcousticArena2);
\draw[<-](TimeFrame) to [bend left=25]node[below,pos=0.5,align=center]{}(AcousticArena);
\draw[->](PerformativeAction) to [bend right=12]node[pos=0.6,above,align=center,xshift=24]{}(TimeFrame);
\draw[->](PerformativeAction2) to [bend left=12]node[pos=0.6,above,align=center,xshift=24]{}(TimeFrame);
\draw[->](MusicalScore)--node[above,pos=0.5]{}(Instruction);
\draw[<-](MusicalScore) to [bend right =12]node[right,pos=0.4,align=center,align=left]{}(PerformativeAction);
\draw[<-](MusicalScore) to [bend left =12]node[right,pos=0.4,align=center,align=left]{}(PerformativeAction2);
\draw[->](PerformativeAction)--node[right,pos=0.5,align=center]{}(AcousticArena);
\draw[->](PerformativeAction2)--node[right,pos=0.5,align=center]{}(AcousticArena2);
\draw[->](Pair)--node[left,pos=0.5,align=center]{}(AcousticArena);
\draw[->](Pair2)--node[left,pos=0.5,align=center]{}(AcousticArena2);
\draw[->](Pair) to [bend right=12]node[left,pos=0.35,align=center]{}(AmbientSounds2);
\draw[->](Pair) to [bend left=12]node[right,pos=0.35,align=center,xshift=5]{}(IncidentalSounds);
\draw[->](Pair2) to [bend right=12]node[left,pos=0.35,align=center]{}(AmbientSounds);
\draw[->](Pair2) to [bend left=12]node[right,pos=0.35,align=center,xshift=5]{}(IncidentalSounds2);
\draw[->](Site) to [bend right=12]node[left,pos=0.5,align=center]{}(AmbientSounds);
\draw[->](Site2) to [bend left=12]node[left,pos=0.5,align=center]{}(AmbientSounds2);
\draw[<-](Listener) to [bend left=12]node[right,pos=0.5,align=center]{}(IncidentalSounds);
\draw[<-](Listener2) to [bend right=12]node[right,pos=0.5,align=center]{}(IncidentalSounds2);
\draw[<-](Listener) to [bend left=20]node[below,pos=0.5,align=center]{}(Site2);
\draw[<-](Listener2) to [bend right=20]node[below,pos=0.5,align=center]{}(Site);
\draw[<-](Sequence2) to [bend right=12]node[right,pos=0.5,align=center]{}(MusicalScoreString);
\draw[<-](Sequence2) to [bend right=12]node[pos=0.4,xshift=20,below,align=left]{}(PerformativeAction);
\draw[<-](Sequence) to [bend left=12]node[pos=0.4,xshift=20,below,align=left]{}(PerformativeAction2);
\draw[<-](MusicalScore) to [bend left=12]node[right,pos=0.6]{}(MusicalScoreString);
\draw[->](PerformativeAction) to [bend left=12]node[left,pos=0.35,align=center]{}(Instruction);
\draw[->](PerformativeAction2) to [bend right=12]node[left,pos=0.35,align=center]{}(Instruction);
\draw[->](Performer) to [bend right=12](MusicalScore);
\draw[->](Performer) to [bend left=12](Sequence);
\draw[->](PerformativeAction) to [bend right=20](MusicalScoreString);
\draw[->](PerformativeAction2) to [bend right=20](Performer);
\end{tikzpicture}
}
\smallskip
\caption[]{Category of elements $\mathcal{A}^{\prime}=\int_{\mathcal{A}}(I)$, for which $\pi_{I}:\mathcal{A}^{\prime}\rightarrow \mathcal{A}$.}
\label{OlogAprime}
\end{figure}
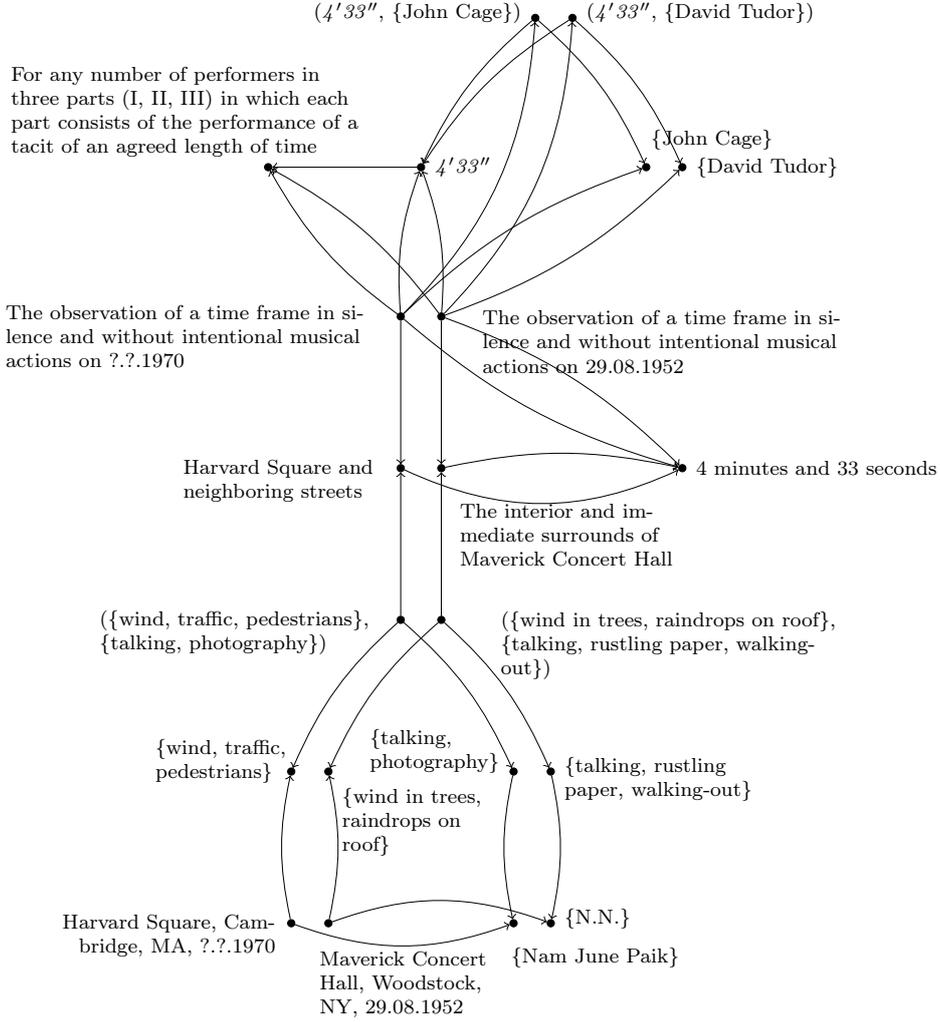

 The discrete opfibration, $(\int_{\mathcal{C}}(I),\pi_{I})$ contains fibers in the form $\pi_{I}^{-1}(s)$, in which $s \in \textrm{Ob}(\mathcal{C})$. This can be seen in Figure~\ref{OlogAprime}, in which a map of instances of \fourmins\space according to the structure of the categorical schema $\mathcal{A}$ aligns instances in $I$ to vertices in $\mathcal{A}$. If the projection functor $\pi_{I}$ is a type of mapping of the category of elements back onto the category $\mathcal{A}$, then the fiber $\pi_{I}^{-1}(s)$ forms a \textit{discrete category}, that is, a category in which there exists only objects and no morphisms between objects. The \textit{pre-image} $\pi_{I}^{-1}(f)$ of $f: s \rightarrow s^{\prime}$ is the set of morphisms between objects $\pi_{I}^{-1}(s)$, and $\pi_{I}^{-1}(s^{\prime})$, for which the subcategory $\pi_{I}^{-1}(f)\subseteq I$ can be described as the function $\pi_{I}^{-1}(f): \pi_{I}^{-1}(s)\rightarrow \pi_{I}^{-1}(s^{\prime})$. An opfibration is \textit{discrete} when there exists a unique morphism emanating from each object in $\pi_{I}^{-1}(s)$. As Spivak asserts, this means that the discrete opfibration contains the same information as $I$ but ``from a different perspective,'' \citep[page 12]{Spivak2013c} which entails the following congruences:
 \begin{align}
 \pi_{I}^{-1}(s)&\cong I(s)\\
 \pi_{I}^{-1}(f)&\cong I(f).
 \end{align}
An opfibration can be viewed as a utility for overlaying the structure of a database back onto a schema. This means that in the case of \fourmins, we can track particular fibers in $\mathcal{A}^{\prime}=\int_{\mathcal{A}}(I)$, given a large enough collection of instances. As an example, consider the fiber $\pi_{I}^{-1}(T)$, where $T$ is `a timestamped performative action.' 
 
 We could image that the function $f: T\rightarrow G$, labelled as the aspect  `temporalizes' is \textit{bijective} in that for each instance of a timestamped performative action, we could expect a unique value for the acoustic arena, even with performances occurring at the same site. While this function suggests that even with the same gesture (performative action), a large set of possible acoustic spaces emerge, the function $t: T\rightarrow Z$, labelled as the aspect `has duration in mins./s.,' is a \textit{non-injective surjective} function. This is given through the fact that for every timestamped performative act, there could be many more \textit{re-enactments} of David Tudor's premiere of the work with the chosen time frame of 4 minutes and 33 seconds.\footnote{As Fetterman notes, there have been many high-profile performances of \fourmins\space as a piano solo of duration four minutes and 33 seconds that also replicate Tudor's gesture from 1952 of opening and closing the piano lid to indicate movements. Such performances have also continued since Cage's death, and have assisted in associating the work only to the time frame of four minutes and 33 seconds, and the instrument of the piano. See \citet[page 80-81]{Fetterman1996}, and \citet[page 290-291]{Labelle2015} for a discussion.} 
 
This suggests that \textit{novel} interpretations of the work are instances that adhere to the division of the work into three sections (the fiber $\pi_{I}^{-1}(T)$ and associated surjective function $u: T\rightarrow D$), but eschew Tudor's timeframe (at object $Z$), while relying on a strong division of sonic materials that can be categorized as incidental ($K$) or ambient ($E$). While the indeterminate nature of the sonic materials present in the work point towards the observation of unique acoustic arenas connected to various physical coverage areas, these areas are inevitably temporalized according to the fiber $\pi_{I}^{-1}(A)$.
 
It follows then that if music is as Cage contends, ``sounds which are just sounds,'' \citep[page 10]{Cage1961} then \fourmins\space demarcates spatiotemporal structures that are \textit{musical}.\footnote{We make this assumption given Brooks' argument that Cage's anechoic experience indicates that ``sound is not suitable to `verify,' [or] to `validate,' the ideas of music or silence; if sound is to be the criterion, music and silence are identical.'' \citep[page 125]{Brooks2007}} This universal condition alerts us to instances of these points in time as sound ecologies that are in principle,  reception mechanisms for articulating acoustic arenas. As Kahn reports, this points to the many performances since Cage's death in 1992 of the ratio between ambient sounds and incidental sounds as a function of the audience reception model the work presents: ``in every performance I've attended the silence has been broken by the audience and become ironically noisy.'' \citep[page 560]{Kahn1997} 
 
\section{Translating $\mathcal{A}$ through two functors}

From Axiom~\ref{Axiom1} we assumed \fourmins, \zeromins, and \one\space are instances of a meta-work Cage called the \textit{Silent piece}. As such, we propose the functors $\phi: \mathcal{A}\rightarrow \mathcal{B}$ and $\psi: \mathcal{A}\rightarrow \mathcal{C}$ as structure preserving maps that will allow us to derive \zeromins, and \one\space from \fourmins. The nomination of the span $\mathcal{C}\xleftarrow{\psi}\mathcal{A}\xrightarrow{\phi}\mathcal{B}$ rather than a strict chronological mapping such as $\mathcal{A}\xrightarrow{\phi}\mathcal{B}\xrightarrow{\kappa}\mathcal{C}$ is a choice informed by Cage's 1990 characterization of \fourmins\space ``becoming'' \one\space\citep[page 577]{Cage2016}, and thus a ``new silent piece'' \citep[page 575]{Cage2016} rather than a direct variation of the previous work \zeromins. 

\subsection{$\phi: \mathcal{A}\rightarrow \mathcal{B}$}

\begin{figure}[t]
\centering
\begin{adjustbox}{width=0.75\textwidth}
\begin{tikzpicture}
\node at (0,0) (MusicalScore) {$\phi(M)$};
\node[right=60mm of MusicalScore](Sequence){$\phi(J)=\phi(L)$};
\node[left =20mm of MusicalScore](Instruction){$\phi(D)$};
\node[below=20mm of MusicalScore](PerformativeAction){$\phi(T)$};  
\node[below=20mm of PerformativeAction](AcousticArena){$\phi(G)$}; 
\node[right=20mm of AcousticArena] (TimeFrame){$\phi(Z)$};
\node[align=left,below=20mm of AcousticArena](Pair){$\phi(A)$};
\node[align=left,below left=20mm and 20mm of Pair](AmbientSounds){$\phi(E)$};
\node[align=left,below right=20mm and 20mm of Pair] (IncidentalSounds){$\phi(K)$};
\node[above=20mm of Sequence] at ($(MusicalScore)!0.5!(Sequence)$) (Performer){$\phi(B)$};
\node[below=15mm of AmbientSounds](Site){$\phi(Q)$};
\draw[thick,<-](TimeFrame)--node[below,pos=0.5]{$\phi(a)$}(AcousticArena);
\draw[thick,->](PerformativeAction)--node[pos=0.6,above,align=center,above,yshift=5]{$\phi(t)$}(TimeFrame);
\draw[thick,->](MusicalScore)--node[above,pos=0.5]{$\phi(c)$}(Instruction);
\draw[thick,<-](MusicalScore)--node[left,pos=0.5,align=center]{$\phi(j)$}(PerformativeAction);
\draw[thick,->](PerformativeAction)--node[left,pos=0.5,align=center]{$\phi(f)$}(AcousticArena);
\draw[thick,->](Pair)--node[right,pos=0.5,align=center]{$\phi(d)$}(AcousticArena);
\draw[thick,->](Pair)--node[left,pos=0.3,align=center,xshift=-5]{$\phi(X)$}(AmbientSounds);
\draw[thick,->](Pair)--node[right,pos=0.3,align=center,xshift=5]{$\phi(Y)$}(IncidentalSounds);
\draw[thick,->](Site)--node[left,pos=0.5,align=center]{$\phi(p)$}(AmbientSounds);
\draw[thick,<-](Sequence.-140) to [bend left=15] node[left,pos=0.5,align=center]{$\phi(h)$}(IncidentalSounds);
\draw[thick,<-,draw](Sequence) -- ++(-90:88mm) to [out=-95,in=0] node[right,pos=0.4,xshift=9]{$\phi(l)$}(Site);
\draw[thick,<-](Sequence)--node[below,pos=0.5,align=center,xshift=5]{$\phi(e)$}(PerformativeAction);
\draw[thick,<-](MusicalScore)--node[left,pos=0.6,xshift=-5]{$\phi(s)$}(Performer);
\draw[thick,<-](Sequence)--node[right,pos=0.6,xshift=5]{$\phi(P)$}(Performer);
\draw[thick,->](PerformativeAction) to [bend left=0]node[right,pos=0.6,xshift=5]{$\phi(w)$}(Performer);
\draw[thick,->](PerformativeAction) to [bend left=0]node[left,pos=0.6,xshift=-5]{$\phi(u)$}(Instruction);
\end{tikzpicture}
\end{adjustbox}
\caption{Category presentation of $\mathcal{B}$ in which the type $\phi(J)$ is given the label `an actant.'}
\label{CategoryB}
\end{figure}
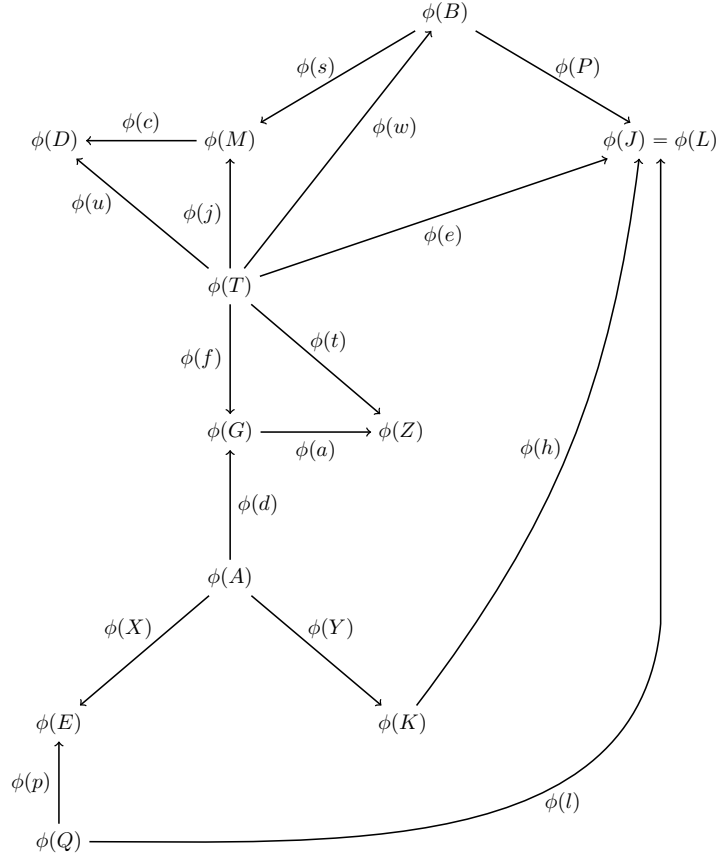

Consider the functor $\phi$, and the categorical schema in Figure~\ref{CategoryB}. As a category presentation of \zeromins\space the functor $\phi$ maps vertices to vertices, arrows to arrows and preserves the path equivalences declared in $\mathcal{A}$ (equations~\ref{eqv1}--\ref{eqv4}). It also maps the vertex $L$ in $\mathcal{A}$ to the vertex $\phi(J)$, such that $\phi(J)=\phi(L)$. 

This mapping also preserves the source and targets of the arrows $\phi(h)$, and $\phi(l)$, and their labelling as `is produced by,' and `contains' respectively from $\mathcal{A}$. We propose the label of the vertex $\phi(J)$ in the olog of $\mathcal{B}$ as the type `an actant.' By actant we mean as in the sense of what \citet{Greimas1983} saw as ``the great functions or roles occupied by the various characters of a narrative.'' \cite[page 505]{Vandendorpe1993} In \zeromins\space these equate to the semiotics of humans in everyday acts. Instances on $\phi(J)$ therefore no longer need be constrained to some \textit{musical qualification} tagged to some person.

Here, the mapping of $\phi(L)$ (`a set of listeners'), and $\phi(J)$ (`a set of performers') to the same vertex implies that the vertex contains \textit{subobjects} as subsets. This can also be described through what \citet[page 61-62]{Spivak2013} calls a \textit{subobject classifier} and \textit{characteristic function}. 

\begin{definition}[Subobject classifier]\label{SubObj}\normalfont
In the category of $\mathbf{Set}$, let the \textit{subobject classifier}, denoted $\Omega$, to be the set $\Omega:=\{\texttt{True},\texttt{False}\}$, together with the monomorphism $\{\ast\}\rightarrow \Omega$ that sends the unique element to \texttt{True}.
\end{definition}
\begin{definition}[Characteristic function]\label{CharFun}\normalfont
Given the subset $S\subseteq C$, we call the corresponding function $p:C\rightarrow \Omega$ the \textit{characteristic function} of $S$ in $C$ such that:
\begin{equation}
p(c) =
\begin{cases}
	\texttt{True}	& \text{if } m(s)=c \textrm{ for some } s \in S \\
	\texttt{False}	& \text{otherwise}
\end{cases}
\end {equation}
\end{definition}

The characteristic function allows for the conceptualization of subobjects of $C$ as being classified by predicates ($p$) on C. This means that for any injective function (monomorphism) between sets $m:S\rightarrowtail C$, there exists a characteristic function on $c$ that returns a value in $\Omega$ (the set of Booleans). If $S$ is a subobject of $C$, then we make $p=\texttt{True}$ if and only if $c$ is found in $S$. We define $C$ to be the vertex in $\mathcal{B}$ labelled $\phi(J)=\phi(L)$ that contains the subset $S$ of $\phi(J)$, and $\phi(L)$ such that:
\begin{equation}
\begin{tikzpicture}
\node at (0,0)(J){$\phi(J)$};
\node at (2,0)(L){$\phi(L)$};
\node at (1,-1.8)(C){$C$};
\node at (3,-1.8)(O){$\Omega$};
\draw[>->](J)--node[left]{$m$}(C);
\draw[>->](L)--node[right]{$m^{\prime}$}(C);
\draw[->](C)--node[below]{$p$}(O);
\end{tikzpicture}
\end{equation}
Here we assert that both a listener and a performer are a type of actant ($m(j)=c$, and $m^{\prime}(l)=c$) and thus $S=\{c\in C\,|\,p(c)=\texttt{True}\}$. Given that some listeners (audience members) assisted in activating the acoustic arena at the premiere of \fourmins\space through generating sounds, we can be classify them in $\mathcal{B}$ as a subset of qualified actants in \zeromins, and thus potential interpreters of the piece. The characteristic function thus also implies that the instance of an instruction in $\phi(D)$ from \zeromins\space becomes \textit{universal} for listeners and performers in $\mathcal{A}$. 

This is because there is a wide range of everyday performative actions that may authenticate the instruction of \zeromins, given the restriction Cage places on the action not being ``the performance of a `musical' composition.'' \citep{Cage1962} Both \citet{Dodd2018} and \citet{Davies1997} contend that this places the work within the realm of performance art rather than music. But Cage notes that, ``indeterminacy, when present in the making of an object . . . is a sign not of identification . . . but simply carelessness with regard to the outcome,'' \citep[page 38]{Cage1961} inferring that music is at least to do with actions and gestures, as it is to do with perceiving \textit{sonic objects}. 

We can think then of the role of the audience in \fourmins\space as an indirect foreshadowing of the role of the performer of in \zeromins---that is, an actant who under particular acoustic conditions, performs some disciplined act whose sounds constitute thematic materials. In the case of \zeromins, these incidental sounds should now be delivered via an amplification system at a maximum setting, for which the act itself should be an  obligation to others rather than a self-centered act. Indeed, in Cage's many documented performances of \zeromins, such obligations included writing correspondences on amplified typewriters, or in the case of the premiere, the writing of the instructions themselves for the work using contact microphones.

The translation $\phi$ also preserves the connection between $L$ and $Q$ (the timestamped site of the performance of the work) through the aspect labelled `contains' ($l$). Though in all documented cases of Cage's performances of \zeromins\space \citep[page 84-90]{Fetterman1996}, these took place in more traditional settings of concert halls or art galleries, as is shown in Cage's 1970 performance of \fourmins\space on Harvard Square, we could expect that the site of \zeromins\space need not be within a traditional concert setting.  

\subsection{$\psi: \mathcal{A}\rightarrow \mathcal{C}$}

\begin{figure}[t]
\centering
\begin{adjustbox}{width=0.81\textwidth}
\begin{tikzpicture}
\node at (0,0) (MusicalScore) {$\psi(M)$};
\node[right=60mm of MusicalScore](Sequence){$\psi(J)=\psi(L)$};
\node[left =20mm of MusicalScore](Instruction){$\psi(D)$};
\node[below=20mm of MusicalScore](PerformativeAction){$\psi(T)$}; 
\node[below=20mm of PerformativeAction](AcousticArena){$\psi(G)$}; 
\node[right=20mm of AcousticArena] (TimeFrame){$\psi(Z)$};
\node[align=left,below=20mm of AcousticArena.west,anchor=west](Pair){$\psi(A)=\psi(E)=\psi(K)$};
\node[above=20mm of Sequence] at ($(MusicalScore)!0.5!(Sequence)$) (Performer){$\psi(B)$};
\node[below=15mm of AmbientSounds,align=left](Site){$\psi(Q)$};
\draw[thick,<-](TimeFrame)--node[below,pos=0.5]{$\psi(a)$}(AcousticArena);
\draw[thick,->](PerformativeAction)--node[pos=0.6,above,align=center,above,yshift=5]{$\psi(t)$}(TimeFrame);
\draw[thick,->](MusicalScore)--node[above,pos=0.5]{$\psi(c)$}(Instruction);
\draw[thick,<-](MusicalScore)--node[left,pos=0.5,align=center]{$\psi(j)$}(PerformativeAction);
\draw[thick,->](PerformativeAction)--node[left,pos=0.5,align=center]{$\psi(f)$}(AcousticArena);
\draw[thick,->](Pair.167.69)--node[right,pos=0.5,align=center]{$\psi(d)$}(AcousticArena);
\draw[thick,<-,draw](Sequence) -- ++(-90:88mm) to [out=-95,in=0] node[right,pos=0.4,xshift=9]{$\psi(l)$}(Site);
\draw[thick,<-](Sequence)--node[below,pos=0.5,align=center,xshift=5]{$\psi(e)$}(PerformativeAction);
\draw[thick,<-](Sequence)--node[right,pos=0.5,align=center,xshift=10]{$\psi(P)$}(Performer);
\draw[thick,<-](MusicalScore)--node[left,pos=0.5,align=center,xshift=-10]{$\psi(s)$}(Performer);
\draw[thick,<-](Sequence.-140) to [bend left=20] node[right,pos=0.6,align=center]{$\psi(h)$}(Pair);
\draw[thick,->](Site) to [bend right=20] node[right,pos=0.5,xshift=5]{$\psi(p)$}(Pair.192.5);
\draw [->,thick] (Pair.north west)++(29:4mm)arc(0:272:0.34)node[left,yshift=18,xshift=-10]{$\psi(\textrm{id}_{A})=\psi(X)=\psi(Y)$};
\draw[thick,->](PerformativeAction)--node[left,pos=0.5,align=center,xshift=-10]{$\psi(u)$}(Instruction);
\draw[thick,->](PerformativeAction)--node[right,pos=0.5,align=center,xshift=10]{$\psi(w)$}(Performer);
\end{tikzpicture}
\end{adjustbox}
\caption{Category presentation of $\mathcal{C}$ in which the morphism $\psi(h)$ is labelled as the aspect `is perceived by,' and the object $\psi(A)$ the type `a set of sounds.'}
\label{CategoryC}
\end{figure}
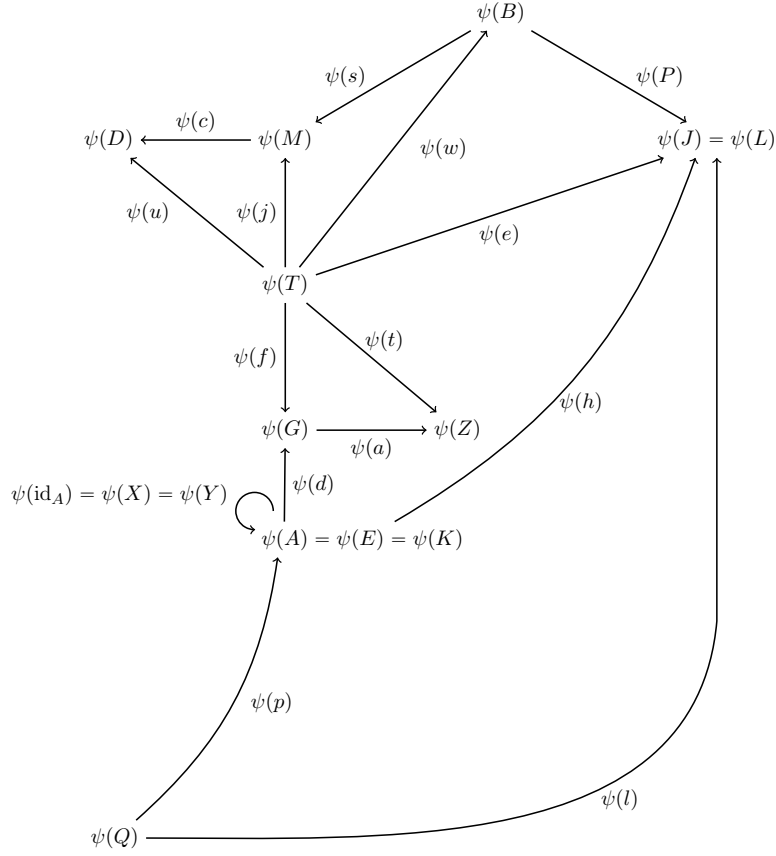

In Figure~\ref{CategoryC} we present the categorical schema $\mathcal{C}$ (\one) as a structure associated to the functor $\psi: \mathcal{A}\rightarrow \mathcal{C}$. As with $\phi$, the functor $\psi$ maps vertices to vertices, arrows to arrows and preserves the path equivalences declared in $\mathcal{A}$ (equations~\ref{eqv1}--\ref{eqv4}). In addition to these path equivalences, given the mapping of objects $E,K\rightarrow \psi(A)$ in $\mathcal{C}$, there is an additional path equivalence of morphisms $\psi(h)\circ \psi(p)\simeq \psi(l)$ as the composition of aspects:
\begin{equation}
(\textrm{is perceived by})\circ(\textrm{produces})\simeq(\textrm{contains}).
\end{equation}
Here, we provide a new label for the aspect $h$ as `is perceived by,' in order to more accurately describe the performance practice Cage suggests in the work. The equivalence asserts that if a timestamped site that produces a set of ambient sounds are perceived by an actant, then that the actant is also located in the site.

The unpublished work \one\space was premiered on November 14, 1989 in Nagoya, and follows \zeromins\space in using amplification in its realization. The manuscript is a fax of a correspondence between Cage and Japanese concert organizers regarding the requirements and instructions for the realization of \one. The instructions state that:

\begin{quote}
You should arrange the sound system so that the whole hall is just on the edge of feedback . . . not actually feeding back but feeling like it might.\\

No piano is necessary. \citep{Cage1989}
\end{quote}

The performative action Cage used in the premiere of \one\space was to direct that the amplification level of the playback system in the concert hall be bought up to near-feedback conditions, and then take a seat in the audience. As in \fourmins, the site of the performance plays heavily in what sounds will be heard, though in \one, the time frame is less well defined, and is an ``unmeasured'' amount of time \citep[page 577]{Cage2016}, which in the premiere was about 12 and a half minutes. The work is thus site-specific given that within the physical arena there may emerge sounds from the local acoustic community, in addition to other ambient sounds found within other remote but intersecting acoustic arenas.

We can think of the mapping of the vertices $E$ (`a set of ambient sounds') and $K$ (`a set of incidental sounds') in $\mathcal{A}$ to $\psi(A)$ in $\mathcal{C}$ as a simple union of tables from $\mathcal{A}$ at $\psi(A)$ (Figure~\ref{CategoryC}). We label the vertex $\psi{A}$, `a set of sounds,' given that we can consider an instance $I(v)$ as something that can be \textit{migrated} between categorical schemas. 

Additionally, given an inherited structure of subobjects that are the $\tau$ projections from the product $A=E\times K$ in $\mathcal{A}$, the sounds that constitute the thematic materials of \one\space are sets of ambient or incidental sounds, or combinations of them. \citet[page 215]{Spivak2013} presents the utility of a \textit{left push forward functor} via a \textit{comma category} in order to describe the union of vertices that as equivalent to parameterized colimits.

\begin{definition}[Comma category]\normalfont
Let $\mathcal{A}$, $\mathcal{B}$ and $\mathcal{C}$ be categories for which there exists the functors $F: \mathcal{A}\rightarrow\mathcal{C}$ and $G:\mathcal{B}\rightarrow\mathcal{C}$. A \emph{comma category} of the functors $F$ and $G$, denoted $(F \downarrow_{\mathcal{C}}G)$ is the category of morphisms of $\mathcal{C}$ from $F$ to $G$ with objects as triples:
\begin{align}
\textrm{Ob}(F\downarrow_{\mathcal{C}}G):=\{(a, b, f)\,|\,a \in \textrm{Ob}(\mathcal{A}) \wedge b \in \textrm{Ob}(\mathcal{B})\wedge\nonumber\\
f:F(a) \rightarrow G(b) \in \textrm{Hom}_{\mathcal{C}}\}.
\end{align}
Let $(a_{1}, b_{1}, f_{1})$ and $(a_{2},b_{2},f_{2})$ be objects of $(F \downarrow_{\mathcal{C}}G)$. The morphisms between these two objects forms the set:
\begin{align}
\{(q,r)\,|\,q:a_{1}\rightarrow a_{2} \in \textrm{Hom}_{\mathcal{A}} \wedge r:b_{1}\rightarrow b_{2} \in \textrm{Hom}_{\mathcal{B}},\nonumber\\
\textrm{such that}\,\,f_{2} \circ F(q)=G(r)\circ f_{1}\}.
\end{align}
for which the following diagram commutes:
\begin{equation}
\begin{tikzpicture}
\node at (0,0) (FAa) {$F(a_{1})$};
\node at (3,0) (FB) {$G(b_{1})$};
\node at (0,-3) (GA) {$F(a_{2})$};
\node at (3,-3) (GB) {$G(b_{2})$};
\draw[->] (FAa)--node[left]{$F(q)$}(GA);
\draw[->] (FB)--node[right]{$G(r)$}(GB);
\draw[->] (FAa)--node[above]{$f_{1}$}(FB);
\draw[->] (GA)--node[below]{$f_{2}$}(GB);
\end{tikzpicture}
\end{equation}
We call the diagram $\mathcal{A}\xrightarrow{F}\mathcal{C}\xleftarrow{G}\mathcal{B}$ the \textit{setup} for the comma category $(F \downarrow_{\mathcal{C}}G)$. There are two \textit{canonical functors} called the \textit{left projection}, denoted $\pi:(F \downarrow_{\mathcal{C}}G)\rightarrow \mathcal{A}$ that sends $(a,b,f)$ to $a$, and a \textit{right projection}, denoted $\tau:(F \downarrow_{\mathcal{C}}G)\rightarrow\mathcal{B}$, that sends $(a,b,f)$ to $b$.
\end{definition}

\begin{definition}[Left pushforward]\normalfont
Let $F:\mathcal{C}\rightarrow\mathcal{D}$ be a functor between categorical schemas, and $I:\mathcal{C}\rightarrow\mathbf{Set}$ an instance on $\mathcal{C}$. We call the functor $\Sigma_{F}(I):\mathcal{D}\rightarrow\mathbf{Set}$ the \textit{left pushforward} on $\mathcal{C}$. To compute $\Sigma_{F}$, we start with an object $d\in\textrm{Ob}(\mathcal{D})$, and form the comma category $(F\downarrow_{\mathcal{D}} d)$ from the setup $\mathcal{C}\xrightarrow{F}\mathcal{D}\xleftarrow{d}\mathbf{1}$, where $\mathbf{1}$ is the terminal category with one object and an identity morphism. There exists a canonical projection functor $\pi:(F\downarrow_{\mathcal{D}} d)\rightarrow \mathcal{C}$, which when composed with $I:\mathcal{C}\rightarrow\mathbf{Set}$ yields the functor $(F\downarrow_{\mathcal{D}} d)\rightarrow\mathbf{Set}$, for which
\begin{equation}
\Sigma_{F}(I)(d):=\colimCC I \circ \pi,
\end{equation}
defines  $\Sigma_{F}(I):\mathcal{D}\rightarrow\mathbf{Set}$ on objects $d\in\mathcal{D}$. Given a morphism in the form $g:d\rightarrow d^{\prime}$, there exists a functor $(F\downarrow_{\mathcal{D}}g):(F\downarrow_{\mathcal{D}}d)\rightarrow(F\downarrow_{\mathcal{D}}d^{\prime})$, such that the following diagram commutes:
\begin{equation}
\begin{tikzpicture}
\node at (0,0)(C0){$(F\downarrow_{\mathcal{D}}d)$};
\node at (5,0)(C1){$(F\downarrow_{\mathcal{D}}d^{\prime})$};
\node at (2.5,-2)(C){$\mathcal{C}$};
\node at (2.5,-3.5)(Set){$\mathbf{Set}$};
\draw[->](C0)--node[above]{$(F\downarrow_{\mathcal{D}}g)$}(C1);
\draw[->](C0)--node[right]{$\pi$}(C);
\draw[->](C1)--node[left]{$\pi^{\prime}$}(C);
\draw[->](C0) to [bend right=20]node[left]{$I\circ \pi$}(Set);
\draw[->](C1) to [bend left=20]node[right]{$I\circ \pi^{\prime}$}(Set);
\draw[->](C)--node[right,pos=0.35]{$I$}(Set);
\end{tikzpicture}
\end{equation}
Given the universal property of colimits we have:
\begin{equation}
\colimCC I \circ \pi \xrightarrow{\Sigma_{F}(I)(g)}\colimCD I \circ \pi^{\prime}.
\end{equation}
\end{definition}

\bigskip
In the case of the vertex $E$ (`a set of ambient sounds'), and the setup $\mathcal{A}\xrightarrow{\psi}\mathcal{C}\xleftarrow{A}\mathbf{1}$, the left push forward $\Sigma_{\psi}(I): \mathcal{C}\rightarrow\mathbf{Set}$ maps $\mathcal{A}$-instances to $\mathcal{C}$-instances. This is achieved through the comma category $(\psi\downarrow_{\mathcal{C}}A)$ with the triple $(E,A,\psi(\textrm{id}_{\mathcal{A}}))$, and a projection ($\pi$) that sends the triple to $A$:
\begin{equation}
\begin{tikzpicture}
\node at (0,0)(C0){$E$};
\node at (5,0)(C1){$A$};
\node at (2.5,-2)(C){$E$};
\node at (2.5,-3.5)(Set){$\mathbf{Set}$};
\draw[->](C0)--node[above]{$\psi(\textrm{id}_{A})$}(C1);
\draw[->](C0)--node[right]{$\pi_{E}$}(C);
\draw[->](C1)--node[left]{$\pi_{A}$}(C);
\draw[->](C0) to [bend right=20]node[left]{$I\circ \pi_{E}$}(Set);
\draw[->](C1) to [bend left=20]node[right]{$I\circ \pi_{A}$}(Set);
\draw[->](C)--node[right,pos=0.35]{$I$}(Set);
\end{tikzpicture}
\end{equation}
which when composed with the instance $I:E\rightarrow \mathbf{Set}$ is a union of instances in the table associated to the vertex $A$ in $\mathcal{C}$ (Figure~\ref{CategoryC}). It is a colimit of the diagram given that for any other object $N$, with morphisms $\xi_{E}: E\rightarrow N$, and $\xi_{A}: A\rightarrow N$, there exists a unique morphism $u: E \rightarrow N$.

\begin{table}\small
\caption{Database table associated to the object $\psi(A)$ in $\mathcal{C}$ as an example of the migration of instances in \fourmins\space to \one}
\label{Table1}
\begin{tabularx}{\textwidth}{| >{\raggedright\arraybackslash}X || >{\raggedright\arraybackslash}X | >{\raggedright\arraybackslash}X |}
\hline
\multicolumn{3}{|c|}{\textbf{A set of sounds}}\\ \hline
\textbf{ID}&\textbf{is perceived by}&\textbf{demarcates}\\\hline
\{wind in trees, raindrops on roof\}&David Tudor&The interior and immediate surrounds of Maverick Concert Hall, NY\\\hline
\{air conditioning, PA system\}&John Cage&The interior of a concert hall in Nagoya Japan\\\hline
\{photography, talking\}&John Cage&Harvard Square and neighboring streets\\\hline
\{talking, rustling paper, walking out\}&N.N.&The interior of Maverick Concert Hall, NY\\\hline
\{talking, rustling paper, walking out, wind in trees, raindrops on roof\}&David Tudor&The interior and immediate surrounds of Maverick Concert Hall, NY\\\hline
\end{tabularx}
\end{table}

Using the same left pushforward on the object $K$ (`a set of incidental sounds'), and $A$ (`a pair of sets $(X,Y)$ in which $X$ are ambient sounds, and $Y$ are incidental sounds') gives:
\begin{align}
\Sigma_{\psi}(I)(\textrm{a set of sounds})=\,&I(\textrm{a set of ambient sounds})\,\sqcup\nonumber\\&I(\textrm{a set of incidental sounds})\,\sqcup\nonumber\\
&I(\textrm{a pair of sets}\, (X,Y)\, \textrm{in which} \,X\, \textrm{are ambient}\nonumber\\ 
&\textrm{sounds, and} \,Y\, \textrm{are incidental sounds}),
\end{align}
for which the associated table of $\psi(A)$ in $\mathcal{C}$ might contain instances such as those found in Table~\ref{Table1}. We provide the new label for the arrow $\psi(h)$ as the aspect `is perceived by' as well as an example $\mathcal{C}$-instance that documents part of the premiere of \one. On one level, the updated labelling of $h$ instantiates the \textit{implied} act of listening found in \fourmins, which in \one\space becomes the principle timestamped performative action. But we have also included in Table~\ref{Table1} the subobject instance N.N., ostensively an audience member from \fourmins. By framing the functor $\psi$ in terms of the migration of instances, we see how the categorical schema $\mathcal{C}$ can describe the role of the listener as an equal actant as the performer, or what Jens Heitjohann calls an ``interdependent'' encounter between audience and performers. \citep{Daniels2016}

\section{Pushout of categorical schemas}\label{PCS}

In order to encapsulate the structure of the thee individual works \one, \zeromins, and \fourmins\space in terms of the \textit{Silent piece} as amalgamation, we propose a construct that takes the two surjective functors $\phi:\mathcal{A}\rightarrow \mathcal{B}$ and $\psi:\mathcal{A}\rightarrow \mathcal{C}$, and derives a \textit{pushout} as colimit through a \textit{free category}. 

\begin{definition}\normalfont
A \emph{pushout} of a diagram, denoted as $S=X\sqcup_{Z}Y$, contains two morphisms: $f : Z \rightarrow X$ and $g : Z \rightarrow Y$ that share a common domain. A pushout is the quotient of $X\sqcup Z\sqcup Y$, or set of equivalence classes:
\begin{equation}
X \sqcup_{Z} Y:= (X \sqcup Z \sqcup Y)/\sim ,
\end{equation} 
such that $\forall z \in Z, z \sim f(z)\,\,\textrm{and}\,\, z \sim g(z)$. In the following communicative square, a $\ulcorner$ denotes the pushout, for which there are the inclusions $i_{1}:Y\rightarrow X \sqcup_{Z}Y$ and $i_{2}:X\rightarrow X \sqcup_{Z}Y$:
\begin{equation}
\begin{tikzpicture}[rotate=270,yscale=-1]
\node at (0,0) (FA) {$X \sqcup_{Z}Y$};
\node[yshift=-15pt,xshift=15pt] at (0,0) {$\ulcorner$};
\node at (3,0) (FB) {$X$};
\node at (0,-3) (GA) {$Y$};
\node at (3,-3) (GB) {$Z$};
\draw[<-] (FA)--node[above]{$i_{1}$}(GA);
\draw[<-] (FB)--node[below]{$f$}(GB);
\draw[<-] (FA)--node[left]{$i_{2}$}(FB);
\draw[<-] (GA)--node[right]{$g$}(GB);
\end{tikzpicture}
\end{equation}
A pushout $S$, and its morphisms $i_{1}, i_{2}$, must also be universal with respect to the communicative square such that for any other triple $(P,j_{1},j_{2})$, there is a unique $u: S\rightarrow P$ called a mediating morphism such that $u \circ i_{2}   = j_{2}$ and $u \circ i_{1}  = j_{1}$:
\begin{equation}\label{pushout}
\begin{tikzpicture}
\node at (0,0)(TL){$S$};
\node[yshift=-15pt,xshift=15pt] at (0,0) {$\ulcorner$};
\node at (3,0)(TR){$Y$};
\node at (0,-3)(BL){$X$};
\node at (3,-3)(BR){$Z$};
\node at (-1.5,1.5)(Q){$P$};
\draw[<-] (TL)--node[above]{$i_{1}$}(TR);
\draw[<-] (TL)--node[left]{$i_{2}$}(BL);
\draw[<-] (BL)--node[below]{$f$}(BR);
\draw[<-] (TR)--node[right]{$g$}(BR);
\draw[<-,dashed] (Q)--node[above]{$u$}(TL);
\draw[<-] (Q) to [bend left=20]node[above]{$j_{1}$}(TR);
\draw[<-] (Q) to [bend right=20]node[left]{$j_{2}$}(BL);
\end{tikzpicture}
\end{equation}
\end{definition}

\begin{definition}[Free category]\normalfont
The \textit{free category} of a directed graph $G$, denoted $Free(G)$, consists of a set of objects that are the vertices, $V$ of $G$, and a set of morphisms that are the identities for each vertex together with finite paths from the set of edges $E$. Composition is defined by concatenation of paths.
\end{definition}

Given the categories $\mathcal{A}$, $\mathcal{B}$, and $\mathcal{C}$, the pushout category of $\phi:\mathcal{A}\rightarrow\mathcal{B}$, and $\psi:\mathcal{A}\rightarrow\mathcal{C}$ contains objects that are generated through an equivalence relation such that $X\in \mathcal{B}$ and $Y\in\mathcal{C}$ are the same image of an object in $\mathcal{A}$. In order to map morphisms, consider the graph $\mathcal{S}^{\prime\prime\prime}$ as
\begin{equation}
\mathcal{S}^{\prime\prime\prime}(X_{1},X_{2})=\bigsqcup_{X_{1}=[Y_{1}],X_{2}=[Y_{2}]}\kern-2em \mathcal{B}(Y_{1},Y_{2})\quad\sqcup\,\bigsqcup_{X_{1}=[Z_{1}],X_{2}=[Z_{2}]} \kern-2em \mathcal{C}(Z_{1},Z_{2}),
\end{equation}
for which $S^{\prime\prime}$ is the free category generated by $\mathcal{S}^{\prime\prime\prime}$, such that a morhpisim $f:X_{1}\rightarrow X_{2}$ is a path $a:X_{1}\rightarrow X_{2}$ in $\mathcal{S}^{\prime\prime\prime}$. We can further impose relations on $S^{\prime\prime}$ to generate $\mathcal{S}^{\prime}$ through functors $\rho:\mathcal{B}\rightarrow S^{\prime}$, and $\kappa:\mathcal{C}\rightarrow S^{\prime}$. This means that the functor identifies a path $(a, b)$ with the composition $(b \circ a)$ when it is defined in $\mathcal{B}$ or $\mathcal{C}$. Similarly, trivial paths as $\textrm{id}_{X}$ are mapped to identity morphisms in $S$ such that $X=[Y]$, where $X$ is in $\mathcal{S}$, and $Y$ is in $\mathcal{C}$, or $\mathcal{B}$. This means that $S^{\prime}$ is the pushout category of $\mathcal{B}$ and $\mathcal{C}$ over objects in $\mathcal{A}$.

In order to check the universal property of pushouts, we apply further relations on $\mathcal{S}^{\prime}$ by asserting the functors $\alpha, \beta:\mathcal{A}\rightrightarrows\mathcal{S}$ must be equal. The functor $F:\mathcal{S}\rightarrow\mathcal{Q}$ is naturally identified with the graph morphism $\gamma: \mathcal{S}^{\prime\prime\prime}\rightarrow U(\mathcal{Q})$, where $U$ is the underlying graph that respects compositions and identity morphisms of $\mathcal{B}$ and $\mathcal{C}$, and coequalizes the two maps from $\mathcal{A}$. This means that the functor $F$ is uniquely identified with a pair of functions from $\mathcal{B}$ and $\mathcal{C}$, which act equally on $\mathcal{A}$. We think of the push out category $\mathcal{S}=\mathcal{B} \sqcup_{\mathcal{A}}\mathcal{C}$ in terms of the following commuting square that describes the relation between \fourmins,\zeromins, and \one\space and meta-work the \textit{Silent piece}:
\begin{equation}
\begin{tikzpicture}[rotate=270,yscale=-1]
\node at (0,0) (FA) {$\mathcal{S}=\mathcal{B} \sqcup_{\mathcal{A}}\mathcal{C}$};
\node[yshift=-15pt,xshift=15pt] at (0,0) {$\ulcorner$};
\node at (3,0) (FB) {$\mathcal{B}$};
\node at (0,-3) (GA) {$\mathcal{C}$};
\node at (3,-3) (GB) {$\mathcal{A}$};
\draw[<-] (FA)--node[above]{$\pi_{1}$}(GA);
\draw[<-] (FB)--node[below]{$\phi$}(GB);
\draw[<-] (FA)--node[left]{$\pi_{2}$}(FB);
\draw[<-] (GA)--node[right]{$\psi$}(GB);
\end{tikzpicture}
\end{equation}

\subsection{Facts as path equivalences in $\mathcal{S}$}

We can also utilize pushouts of morphisms between objects in $\mathcal{S}$ in terms of what \citet[page 11]{Spivak2012} nominate as \textit{facts}, or the semantics of path equivalences. We have already encountered the semantics of path equivalences regarding the composition of morphisms, but they are also useful when reasoning on schemas with pushouts. Consider the olog of the categorical schema $\mathcal{S}$ of the meta-work the \textit{Silent piece} given in Figure~\ref{OlogS}.

\begin{figure}[t!]
\begin{adjustbox}{width=\textwidth}
\begin{tikzpicture}
	\tikzset{
		A-node/.style={
		draw=black,draw,inner sep=3pt,minimum 												size=8mm,label={[align=left]90:#1}},
		}
\node[A-node,label={[align=left]100:$M_{\mathcal{C}}=M_{\mathcal{B}}$},align=left] at (0,0) (MusicalScore) {a musical\\score};
\node[A-node,label={[align=left]80:$P_{\mathcal{C}}=P_{\mathcal{B}}$},right=60mm of MusicalScore.east,align=left](Sequence){a set of actants};
\node[A-node,label={[align=left]90:$D_{\mathcal{C}}=D_{\mathcal{B}}$},left =20mm of MusicalScore](Instruction){an instruction};
\node[A-node, label={[align=left,xshift=-5,yshift=8]30:$T_{\mathcal{C}}=T_{\mathcal{B}}$},below=20mm of MusicalScore,align=left](PerformativeAction){a timestamped\\performative action}; 
\node[A-node={$\hspace{-25mm}G_{\mathcal{C}}=G_{\mathcal{B}}$},below=20mm of PerformativeAction](AcousticArena){an acoustic arena}; 
\node[A-node={$Z_{\mathcal{C}}=Z_{\mathcal{B}}$},right=15mm of AcousticArena,align=left,fill=white] (TimeFrame){a time\\frame};
\node[A-node={$\hspace{-25mm}A_{\mathcal{C}}=A_{\mathcal{B}}$},align=left,below=25mm of AcousticArena](Pair){a set of sounds};
\node[A-node,label={[align=left]70:$E_{\mathcal{B}}$},align=left,below left=20mm and 5mm of Pair](AmbientSounds){a set of  ambient\\sounds};
\node[A-node={$K_{\mathcal{B}}$},align=left,below right=20mm and 5mm of Pair] (IncidentalSounds){a set of incidental\\sounds};
\node[A-node={$\hspace{-37mm}Q_{\mathcal{C}}=Q_{\mathcal{B}}$},below=20mm of AmbientSounds,align=left](Site){a timestamped\\site};
\node[A-node,above=35mm of PerformativeAction,label={[align=left]90:$B_{\mathcal{C}}=B_{\mathcal{B}}$},align=left,text width=40mm] at ($(MusicalScore)!0.5!(Sequence)$) (Performer){a pair $(s,P)$, where $s$ is a musical score and $P$ is a set of actants who interpret $s$};
\begin{pgfonlayer}{background}
\draw[thick,<-](TimeFrame)--node[below,pos=0.5,align=center]{\small is\\\small actualized\\\small during}(AcousticArena);
\draw[thick,->](PerformativeAction)--node[pos=0.65,above,align=center,above,xshift=20]{\small has a duration\\\small in mins./s.}(TimeFrame);
\draw[thick,->](MusicalScore)--node[above,pos=0.5]{\small contains}(Instruction);
\draw[thick,<-](MusicalScore)--node[left,pos=0.4,align=right]{\small is realization\\\small of}(PerformativeAction);
\draw[thick,->](PerformativeAction)--node[left,pos=0.4,align=center]{\small temporalizes}(AcousticArena);
\draw[thick,->](Pair)--node[right,pos=0.4,align=center]{\small demarcates}(AcousticArena);
\draw[thick,->](Pair)--node[left,pos=0.3,align=center,xshift=-5]{\small{includes}}(AmbientSounds.30);
\draw[thick,->](Pair)--node[right,pos=0.3,align=center,xshift=5]{\small{includes}}(IncidentalSounds);
\draw[thick,->](Site)--node[left,pos=0.7,align=center]{\small produces}(AmbientSounds);
\draw[thick,<-](Sequence.248) to [bend left=15] node[right,pos=0.8,align=left,xshift=5]{\small is produced\\\small by}(IncidentalSounds);
\coordinate[xshift=0,yshift=0] (P1) at (IncidentalSounds-|Sequence);
\draw[thick,<-,draw](Sequence) -- ++(-90:110mm) to [out=-95,in=0] node[right,pos=0.4,xshift=9]{\small contains}(Site);
\draw[thick,<-](Sequence)--node[right,pos=0.55,align=center,yshift=-15,align=left]{\small is enacted\\\small by}(PerformativeAction);
\draw[thick,<-](Sequence.227) to [bend left=27]node[left,pos=0.3,xshift=-5,align=right]{\small is perceived\\\small by}(Pair);
\draw[thick,->](Site) to [bend right=20] node[right,pos=0.3,xshift=5]{\small thematizes}(Pair);
\draw[thick,->](PerformativeAction)--node[right,pos=0.45,align=left,yshift=17,text width=25mm,xshift=33]{\small yields, insofar as it is a realization of a musical score and enacted by a set of actants}(Performer.south);
\draw[thick,->](Performer)--node[right,pos=0.4,align=center,xshift=10]{\small yields as $P$}(Sequence.north);
\draw[thick,->](Performer)--node[left,pos=0.4,align=center,xshift=-10]{\small yields as $s$}(MusicalScore.north);
\draw[thick,->](PerformativeAction)--node[left,pos=0.35,align=center,xshift=-9]{\small fulfills}(Instruction);
\draw[thick,->](AcousticArena) to [bend right=20]node[left,pos=0.35,align=center,xshift=-5,yshift=-6]{\small contains}(AmbientSounds.110);
\node[align=left,A-node,{label={[align=left]100:$W$}}] at (AcousticArena-|Instruction) (W){a localized\\acoustic typology};
\node[align=left,A-node,{label={[align=left]100:$L$}},left=13mm of W,yshift=0] at (Instruction|-Pair) (L){an enfolded\\acoustic field};
\draw[thick,->](Instruction)--node[left,pos=0.5,xshift=-2]{determines}(W);
\draw[thick,->,font=\small](AcousticArena)--node[below,pos=0.5,yshift=-15]{instantiates}(W);
\draw[thick,->,align=right,font=\small](W) to [bend right=0]node[left,pos=0.5,xshift=-8,yshift=5]{is\\embedded\\in}(L);
\draw[thick,->,align=left,font=\small](AmbientSounds)--node[right,pos=0.5,yshift=19,xshift=-4]{is\\propagated\\through}(L);
\draw[thick,->,font=\small](Site) to [bend left=20]node[left,pos=0.5,xshift=-5]{spatializes}(L);
\end{pgfonlayer}
\end{tikzpicture}
\end{adjustbox}
\caption[]{Olog of the pushout category $\mathcal{S}$ with additional chained pushouts $W=G\bigsqcup_{T}D$, and $L=W\bigsqcup_{D}T$, and morphism `spatializes' ($m: Q\rightarrow L$).}
\label{OlogS}
\end{figure}
We present the pushouts of $W=G\bigsqcup_{T}D$, and $L=W\bigsqcup_{D}T$ as a means in which to reason on how `a timestamped performative action' ($T$) maps to both `an instruction' ($D$) and `an acoustic arena' ($G$), and secondly, the functions from $G$ to `a set of ambient sounds' ($E$) and to the pushout $W$, the morphism `instantiates.' We label $W$ as `a localized acoustic typology,' 
by which we mean in the sense of a \textit{building typology} in the field of architecture: that is, a set of acoustic spaces that are united through some structural, functional or formal characteristics. In the case of $W$, these fundamental traits emerge from the equivalence relation between the functions $f$ (`temporalizes'), and $u$ (`fulfills') in which a (musical) instruction is paired to an acoustic arena in $W$. 

In the case of $L$, we can think of how these particular typologies of the \textit{Silent piece} are both localized in the time domain as well as spatialized within a larger structure. As such, we give $L$ the label `an enfolded acoustic field.' Here, we consider that the presence of ambient sounds ($E$) that emerge from a site at a point in time ($Q$), form a background acoustic condition that \textit{enfolds} the emerging typologies. This mimics the classic architectonic condition of an \textit{interior} versus an \textit{exterior}. \citep{Fowler2012} 

This implies that we have a localization of acoustic typologies in the time domain, and an embedding of the typologies within an acoustic field bound to a spatial domain. This is specifically expressed in $\mathcal{S}$ in the additional morphism `spatializes,' that connects `a timestamped site' to `an enfolded acoustic field,' and is equivalent to the composition: (is propagated through)$\,\circ\,$(produces). 

But there is also a significance to the chaining of pushouts $W$ and $L$. This can be found in the fact that given $W$ as the pushout of the morphisms `fulfills' and `temporalizes,' and $L$ as the pushout of the morphisms `contains' and `determines,' it follows that $L$ is also the pushout of `temporalizes' then `contains' and `fulfills.'\footnote{See \citet[page 132]{Schubert1972} for proof via the dual \textit{pasting law} for pullbacks.}


\section{Fiber order}
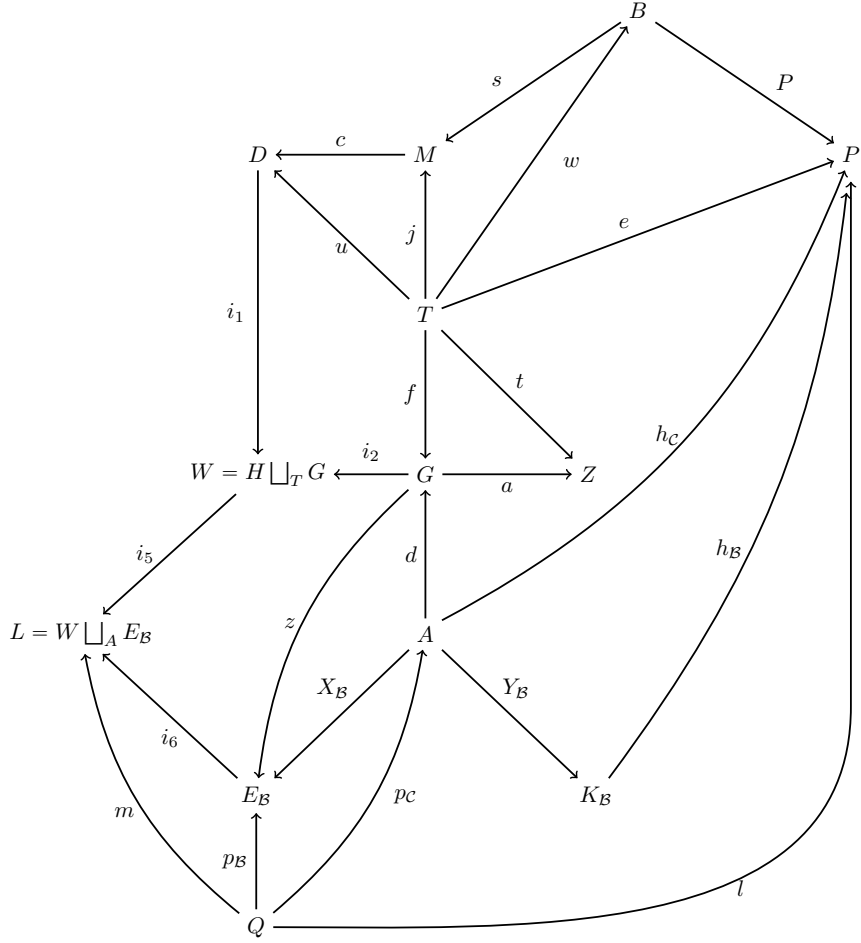
\begin{figure}[t!]
\begin{adjustbox}{width=0.9\textwidth}
\begin{tikzpicture}[cross line/.style={preaction={draw=white, -,line width=4pt}}]
\small
\node at (0,0) (MusicalScore) {$M$};
\node[right=60mm of MusicalScore](Sequence){$P$};
\node[left =20mm of MusicalScore](Instruction){$D$};
\node[below=20mm of MusicalScore](PerformativeAction){$T$}; 
\node[below=20mm of PerformativeAction](AcousticArena){$G$}; 
\node[right=20mm of AcousticArena] (TimeFrame){$Z$};
\node[align=left,below=20mm of AcousticArena](Pair){$A$};
\node[align=left,below left=20mm and 20mm of Pair](AmbientSounds){$E_{\mathcal{B}}$};
\node[align=left,below right=20mm and 20mm of Pair] (IncidentalSounds){$K_{\mathcal{B}}$};
\node[below=15mm of AmbientSounds](Site){$Q$};
\node[above=20mm of MusicalScore] at ($(MusicalScore)!0.5!(Sequence)$)(Performer){$B$};
\node at (AcousticArena-|Instruction) (W){$W=H\bigsqcup_{T}G$};
\node[left=15mm of W] at (Instruction|-Pair) (L){$L=W\bigsqcup_{A}E_{\mathcal{B}}$};
\draw[thick,<-](TimeFrame)--node[below,pos=0.5]{$a$}(AcousticArena);
\draw[thick,->](PerformativeAction)--node[pos=0.6,above,align=center,above,yshift=5]{$t$}(TimeFrame);
\draw[thick,->](MusicalScore)--node[above,pos=0.5]{$c$}(Instruction);
\draw[thick,<-](MusicalScore)--node[left,pos=0.5,align=center]{$j$}(PerformativeAction);
\draw[thick,->](PerformativeAction)--node[left,pos=0.5,align=center]{$f$}(AcousticArena);
\draw[thick,->](Pair)--node[left,pos=0.5,align=center]{$d$}(AcousticArena);
\draw[thick,->](Pair)--node[left,pos=0.3,align=center,xshift=-5]{$X_{\mathcal{B}}$}(AmbientSounds);
\draw[thick,->](Pair)--node[right,pos=0.3,align=center,xshift=5]{$Y_{\mathcal{B}}$}(IncidentalSounds);
\draw[thick,->](Site)--node[left,pos=0.5,align=center]{$p_{\mathcal{B}}$}(AmbientSounds);
\draw[thick,<-,shorten <=10](Sequence) to [bend left=15] node[left,pos=0.6,align=center]{$h_{\mathcal{B}}$}(IncidentalSounds);
\draw[thick,<-,draw,shorten <=5](Sequence) -- ++(-90:86mm) to [out=-90,in=0] node[right,xshift=5,pos=0.4]{$l$}(Site);
\draw[thick,<-](Sequence)--node[left,pos=0.5,align=center,yshift=5]{$e$}(PerformativeAction);
\draw[thick,<-](Sequence) to [bend left=20]node[left,pos=0.5,xshift=0]{$h_{\mathcal{C}}$}(Pair.40);
\draw[thick,->](Site) to [bend right=20] node[right,pos=0.5,xshift=5]{$p_{\mathcal{C}}$}(Pair);
\draw[thick,->](PerformativeAction)--node[below,pos=0.5]{$u$}(Instruction);
\draw[thick,->](Instruction)--node[left,pos=0.5,xshift=-2]{$i_{1}$}(W);
\draw[thick,->](AcousticArena)--node[above,pos=0.5,yshift=2]{$i_{2}$}(W);
\draw[thick,->](AcousticArena) to [bend right=20]node[left,pos=0.5,yshift=-2]{$z$}(AmbientSounds);
\draw[thick,->](W)--node[left,pos=0.5,xshift=-3]{$i_{5}$}(L);
\draw[thick,->](AmbientSounds)--node[below,pos=0.5,yshift=-2]{$i_{6}$}(L);
\draw[thick,->](Site) to [bend left=20]node[below,pos=0.5,xshift=-5]{$m$}(L);
\draw[thick,<-](Sequence)--node[right,pos=0.5,align=center,xshift=10]{$P$}(Performer);
\draw[thick,<-](MusicalScore)--node[left,pos=0.5,align=center,xshift=-10]{$s$}(Performer);
\draw[thick,->](PerformativeAction)--node[right,pos=0.5,align=center,xshift=10]{$w$}(Performer);
\end{tikzpicture}
\end{adjustbox}
\caption{Category presentation $\mathcal{S}$.}
\label{OlogSprime}
\end{figure}

In order to extract a semantics of the spatio-temporal order inherent in $\mathcal{S}$, we can further investigate what \citet[page 20]{Spivak2012} call the \textit{language} of an olog: that is, its \textit{fiber order}. An olog's language consists of a set of path equivalences (as individual facts $\epsilon$) from the olog's underlying graph $G$, which is denoted $\mathbf{eqn}(G)$. An example from category $\mathcal{S}$ (Figure~\ref{OlogSprime}) is the equation $u\simeq c \circ j$, which asserts that the realization of a musical score consists of a fulfillment of the instruction of the score (cf. Figure~\ref{OlogS}). 

We denote the set $\mathbf{eqn}(G)$ as a collection of all possible equivalence statements, for which we call a subset of these equations the $G$-specification, $E\subseteq \mathbf{eqn}(G)$. The specification on the category $\mathcal{S}$ is then given as $\mathbf{spec}(\mathcal{S})\langle G,E\rangle$, which is the collection of all facts that hold in $\mathcal{S}$. For the \textit{Silent piece}, we nominate the following collection of facts as comprising $E$:
\begin{align}\label{E}
&E_{1}:= u \simeq (c \circ j) \simeq (c\circ s \circ w) \nonumber\\
&E_{2}:= e \simeq (P \circ w) \nonumber\\
&E_{3}:= j \simeq (s\circ w)\nonumber\\
&E_{4}:= t \simeq (a \circ f)\nonumber\\
&E_{5}:= (i_{1} \circ u) \simeq (i_{2} \circ f) \simeq (i_{1} \circ c \circ j) \simeq (i_{1} \circ c \circ s \circ w)\nonumber\\
&E_{6}:=(i_{5} \circ i_{1} \circ u) \simeq (i_{5} \circ i_{1} \circ c \circ j) \simeq (i_{5} \circ i_{1} \circ c \circ s \circ w) \simeq\nonumber\\
&\qquad\qquad(i_{5} \circ i_{2} \circ f) \simeq (i_{6} \circ z \circ f)\nonumber\\
&E_{7}:= (i_{5} \circ i_{2} \circ d) \simeq (i_{6} \circ X_{\mathcal{B}}) \simeq (i_{6} \circ z \circ d)\nonumber\\
&E_{8}:= h_{\mathcal{C}} \simeq (h_{\mathcal{B}} \circ Y_{\mathcal{B}})\nonumber\\
&E_{9}:= m \simeq (i_{6}\circ p_{\mathcal{B}}) \simeq (i_{6} \circ X_{\mathcal{B}} \circ p_{\mathcal{C}}) \simeq (i_{6} \circ z \circ d \circ p_{\mathcal{C}}) \simeq \nonumber\\
&\qquad\qquad (i_{5}\circ i_{2} \circ d \circ p_{\mathcal{C}})\nonumber\\
&E_{10}:= l \simeq (h_{\mathcal{C}} \circ p_{\mathcal{C}}) \simeq (h_{\mathcal{B}} \circ Y_{\mathcal{B}} \circ p_{\mathcal{C}})\nonumber\\
&E_{11}:=p_{\mathcal{B}} \simeq (X_{\mathcal{B}} \circ p_{\mathcal{C}}) \simeq (z \circ d \circ p_{\mathcal{C}})\nonumber\\
&E_{12}:= X_{\mathcal{B}} \simeq (z \circ d)\nonumber \\
&E_{13}:= (i_{5} \circ i_{2}) \simeq (i_{6} \circ z)\nonumber \\
\end{align}

From $E$ we can also consider how particular instances that form the category of elements $\int_{\mathcal{S}}(I)$ \textit{satisfy} (model) an $\mathcal{S}$-fact $\epsilon \in \mathbf{eqn}(\mathcal{S})$ as $I\models_{\mathcal{S}} \epsilon$. Consider the population of $\mathcal{S}$ depicted in Figure~\ref{OlogSprime} with instances from any of the three premieres of the works \fourmins, \zeromins, and \one, such that there exists a closure operator, $(-)^{\bullet}$, denoted $E^{\bullet}$ of the $\mathcal{S}$-specification, which is the set of all entailed equations. Thus, given $I\models_{\mathcal{S}} E$, such that $E$ entails a $G$-fact $\epsilon$, denoted $E\vdash_{\mathcal{S}} \epsilon$, any model of the specification then satisfies the fact. As an example, this means that in the case of $E_{12}^{\bullet}=X_{\mathcal{B}} \vdash (z \circ d)$, we can assert for any instance in $\int_{\mathcal{S}}(I)$ of a set of sounds (thematic materials), there exists a subset of ambient sounds that are contained within an acoustic arena that is itself demarcated through the parent set of sounds.



\begin{figure}
\centering
\adjustbox{max width=\textwidth}{
\begin{tikzpicture}
\node at (0,0) (E6) {$E_{6}$}; 
\node at (0mm, -14mm) (E5) {$E_{5}$}; 
\node at (-20mm, -28mm) (E7) {$E_{7}$}; 
\node at (-28mm, -14mm) (E9) {$E_{9}$}; 
\node at (-42mm, -14mm) (E11) {$E_{11}$}; 
\node at (0mm, -28mm) (E1) {$E_{1}$}; 
\node at (-14mm, -42mm) (E13) {$E_{13}$}; 
\node at (14mm, -28mm) (E10) {$E_{10}$}; 
\node at (0mm, -42mm) (E3) {$E_{3}$}; 
\node at (28mm, -42mm) (E4) {$E_{4}$}; 
\node at (42mm, -42mm) (E2) {$E_{2}$}; 
\node at (14mm, -42mm) (E8) {$E_{8}$}; 
\node at (0mm, -60mm) (emptyset) {$\varnothing$};
\node at (E11|-E13) (E12) {$E_{12}$}; 
\node at (0mm, 16mm) (top) {$E$};
\begin{pgfonlayer}{background}
\draw[thick](top.south)--(E6);
\draw[thick](top.south)--(E9.north);
\draw[thick](top.south)--(E11.north);
\draw[thick](top.south)--(E10);
\draw[thick](top.south)--(E2);
\draw[thick](top.south)--(E4);
\draw[thick](E6.south)--(E13.north);
\draw[thick](E11.south)--(E12.north);
\draw[thick](E9.south)--(E12.north);
\draw[thick](E7.north)--(E9.south);
\draw[thick](E7.south)--(E13.north);
\draw[thick](E13)--(emptyset.north);
\draw[thick](E6.south)--(E5);
\draw[thick](E5)--(E1);
\draw[thick](E1)--(E3);
\draw[thick](E10)--(E8);
\draw[thick](E12.north)--(E7.south);
\draw[thick](emptyset.north)--(E12.south);
\draw[thick](emptyset.north)--(E3);
\draw[thick](emptyset.north)--(E4);
\draw[thick](emptyset.north)--(E2.south);
\draw[thick](emptyset.north)--(E8.south);
\end{pgfonlayer}
\end{tikzpicture}
}
\caption[]{Lattice of the fiber order, $\mathbf{fbr}(\mathcal{S})=\langle\mathbf{spec}(\mathcal{S}),\geq_{\mathcal{S}}\rangle$ of $E$.}
\label{Lattice1}
\end{figure}
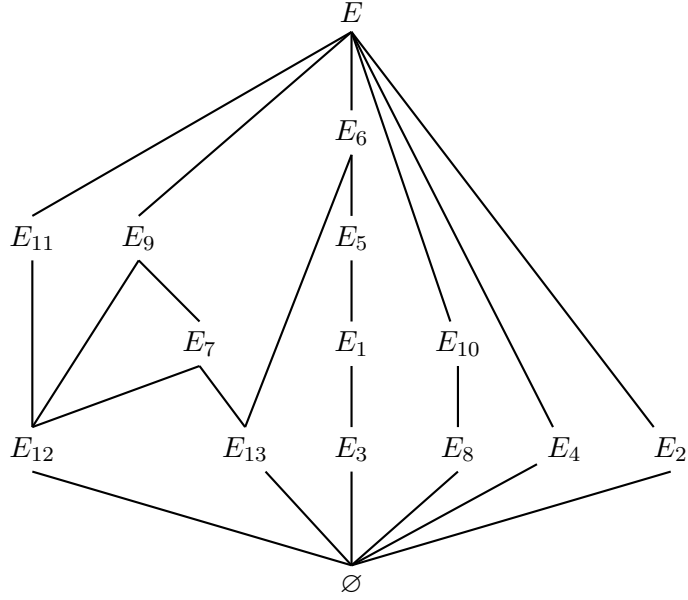

We can further consider an entailment oder of subsets of $E$ as a lattice (Figure~\ref{Lattice1}). The lattice describes the fiber order, denoted $\mathbf{fbr}(\mathcal{S})=\langle\mathbf{spec}(\mathcal{S}),\geq_{\mathcal{S}}\rangle$, and visualizes the hierarchy of facts found in equation~\ref{E}. As an example, we say that $E_{3}\leq_{\mathcal{S}}E_{1}$, given the equations in $E_{1}$ entail the equations in $E_{3}$, or duly that $E_{1}^{\bullet}\supseteq E_{3}$. The fiber order of $\mathcal{S}$ also allows us to understand how particular fibers of the \textit{Silent piece} entail particular acoustic consequences that can be contextualized through spatio-temporal relations. 

Consider the fact $E_{6}$ (cf. equation~\ref{E}), such that $E_{6}^{\bullet}\supseteq (E_{5}, E_{1},E_{3},E_{13})$, and thus all equations entailed in $E_{6}$ cover those found in $E_{5}, E_{1},E_{3},E_{13}$ (Figure~\ref{Lattice1}). Semantically, we can then assert that given that $E_{6}$ is the collection of fibers between $T=$ `timestamped performative action,' and $W=$ `an enfolded acoustic field,' it also entails that:
\begin{itemize}
\item a timestamped performative action fulfills an instruction that determines a localized acoustic typology, while also temporalizing an acoustic arena that is realized by an actant interpreting a musical score ($E_{5}, E_{1}, E_{3}$).
\item an acoustic arena that contains a set of ambient sounds that propagate through an enfolded acoustic field is the instantiation of a localized acoustic typology that is embedded within the enfolded acoustic field ($E_{13}$). 
\end{itemize}

From the fact $E_{9}$, whose paths in $\mathcal{S}$ link $Q=$ `a timestamped site' to $L=$ `an enfolded acoustic field,' we have the relation of $E_{9}^{\bullet}\supseteq (E_{12},E_{7},E_{13})$. Here we ascertain that the timestamped site spatializes an enfolded acoustic field through its production of ambient sounds and subsequent thematization of all sounds present. But further to this, we also have: 
\begin{itemize}
\item a set of sounds containing an ambient subset demarcates an acoustic arena which instantiates a local acoustic typology that is embedded within an enfolded acoustic field, for which the ambient subset of sounds similarly propagate through the enfolded field ($E_{7},E_{12}$). 
\end{itemize}

Finally, let us consider $E_{8}\leq_{\mathcal{S}}E_{10}$, and the relationship between $Q=$ `a timestamped site,' and the containment (the aspect $l$) of $P=$ `a set of actants' within the site. From $E_{10}^{\bullet}\supseteq E_{8}^{\bullet}$ we have $(h_{\mathcal{C}} \circ p_{\mathcal{C}})\vdash h_{\mathcal{C}}$, and $(h_{\mathcal{B}} \circ Y_{\mathcal{B}} \circ p_{\mathcal{C}})\vdash (h_{\mathcal{B}} \circ Y_{\mathcal{B}})$, such that:
\begin{itemize}
\item{a timestamped site thematizes a set of sounds that are perceived by a set of actants for which there exists a subset of sounds that are incidental, and are produced by a set of actants}
\end{itemize}

\section{Discussion}

As is evident from the hierarchy of $\mathbf{fbr}(\mathcal{S})$ in Figure~\ref{Lattice1}, the facts $E_{6}, E_{9}$ and $E_{10}$ cover a number of sub-facts, and describe those most persistent structures of the \textit{Silent piece}. The semantics of these facts generalize how architectonic articulations, program, landscape and agency are facets of a spatio-temporal system.

In terms of the interior/exterior condition we see what \citet{Joseph1997} identifies as Cage's conceptualization of the \textit{interpenetration} of inside and outside, and the influence of Moholy-Nagy's (\citeyear{Moholy-Nagy2012}) definition of architectural space. 
The spatio-temporal dynamics of the categorical schema $\mathcal{S}$ describe an interior as not singularly mapped to any particular architectonic articulation of space. Here, an acoustic arena, together with an instruction are a pair that map to a acoustic typology, such that local spatial articulations of the site of performance need not be uniform across different performances with regard to materiality, volume or containment for there to exist an acoustic typology that unites them. 

Cage's notion of interpenetration is also captured through the pushout $W \in \mathcal{S}$ (`a localized acoustic typology') as a colimit, and thus an object obtained by sewing together $D$ (`an instruction') and $G$ (`an acoustic arena') according to the morphisms `determines,' and `instantiates.' The fact that the second pushout $L$ (`an enfolded acoustic field') sews the objects $W$ and $E_{\mathcal{B}}$ (`a set of ambient sounds') via the morphisms `instantiates,' and `contains,' means that $L$ covers not only an acoustic arena, but similarly a performative action. We can therefore think of an interior as a \textit{typing} of acoustic arenas across time through a particular performative action, for which the space in which these types converge is an enfolded acoustic field as exterior. 

Regarding the role of actants, there has been much discussion in the literature concerning the true function of auditor and performer in \fourmins\space\citep{Labelle2015,Davies2003,Kim-Cohen2009,Kotz2009}. We showed that Spivak's left pushforward functor $\Sigma_{F}(I):\mathcal{D}\rightarrow\mathbf{Set}$ on a category $\mathcal{C}$ where $F: \mathcal{C}\rightarrow\mathcal{D}$, and $I:\mathcal{C}\rightarrow\mathbf{Set}$, is an operation that allows the migration of instances between schemas (such as an audience member of \fourmins\space mapping to a performer of \zeromins). This suggests the roles of actants as interchangeable within $\mathcal{S}$---which is exemplified in the subobject classifier and characteristic function $p: P\rightarrow \Omega$. We see this also expressed in $\mathcal{S}$ through the fact that timestamped performative actions may map to a set of actants ($e:T\rightarrow P$) who perceive a set of sounds ($h_{\mathcal{C}}:A\rightarrow P$), yet are not intentional producers of a subset of these sounds that are characterized as incidental (the case in \fourmins\space and \one). 

It may be equally true that actants that are agents of a performative action may both perceive a set of sounds as well as produce a subset that are incidental (that is, intentional as in \zeromins). But as Cage's anechoic experience suggests, in all cases actants will be \textit{unintentional} producers of incidental sounds (however subtle) such as those produced by the human body. This condition further highlights Brooks' notion of Cage's \textit{pragmatics of silence}, in which there does not exist an environment in which sounds are not heard or produced, and therefore ``sound is not suitable to `verify,' [or] to `validate,' the ideas of music or silence.'' \citep[page 125]{Brooks2007} 

The final persistent structure of $\mathbf{fbr}(\mathcal{S})$ is the notion of site in the \textit{Silent piece}, and how it emerges not simply through a localization of the performative act according to a shared timestamp, but through the functions of spatialization, production, containment, and thematization. When we consider the full subcategory $\mathcal{S}^{\ast}$ as consisting of the objects $(L,E_{\mathcal{B}},A,P,Q)$ from $\mathcal{S}$, then the object $Q$ (`a timestamped site') acts as the limit of this diagram. This points towards how the presence of $Q$ in $\mathcal{S}^{\ast}$ becomes a type of optimal solution for which each of the objects in the diagram is parametrized through the morphisms: `spatializes,' `produces,' `thematizes,' and `contains.'

This aligns to what \citet[page 13]{Kruse2019} and \citet[page 9]{Massey2005} identify as the multiplicity and plurality of space, and how a site is something that is constantly being generated through interactions. In the case of the \textit{Silent piece}, the site is both a container to actants, while at the same time assigns a set of sounds a thematic function as musical materials, and spatializes a subset of these sounds that are indigenous in an enfolded acoustic field.

\section{Conclusion}

We have drawn on the work of Spivak in applied category theory in order to provide a mathematical context to the question of how to derive a semantics of Cage's meta-work the \textit{Silent piece}. The database approach to representing instances from previous premieres through a category of instances led to considering a fiber order on a push out category whose formal language provides a way in which to reason on Cage's \textit{pragmatics of silence} (Axiom~\ref{Axiom1}).

The semantics of the \textit{Silent piece} that we have derived points towards an \textit{auditory awareness} and resultant audience reception model that identifies an underlying and therefore embedded ecology of sound within the landscape. It further highlights how Cage instrumentalizes this awareness through the declaration that the agency of listening, is concurrent to the agency of sound production---whether intentional or unintentional. When Cage asserts that silence is simply ``the multiplicity of activity that constantly surrounds us,'' \citep[page 245]{KostelanetzCage2003} he is also framing how the meta-work of the \textit{Silent piece} maps actants to an attentiveness of a site regardless of its spatially rarefied qualities, or the specific nature of the performative act being undertaken.

We consider the utility of Spivak's ontological logs within the context of database schemas as a valuable method for analyzing Cage's indeterminate structures: particularly in regard to scores that are communicated solely through natural language instructions. We anticipate that the mathematical context of categorical schemas may also open up other research programs for investigating both philosophical as well as compositional aspects of a variety of musical systems. 



\appendix

\bibliography{Database}

\begin{thebibliography}{46}
\providecommand{\natexlab}[1]{#1}
\providecommand{\url}[1]{\texttt{#1}}
\expandafter\ifx\csname urlstyle\endcsname\relax
  \providecommand{\doi}[1]{doi: #1}\else
  \providecommand{\doi}{doi: \begingroup \urlstyle{rm}\Url}\fi

\bibitem[Awodey(2010)]{Awodey2010}
S.~Awodey.
\newblock \emph{Category Theory}.
\newblock Oxford Logic Guides. OUP Oxford, Oxford, 2010.

\bibitem[Blesser and Salter(2007)]{BlesserSalter2007}
Barry Blesser and Linda-Ruth Salter.
\newblock \emph{Spaces speak are you listening? Experiencing Aural
  Architecture.}
\newblock MIT Press, Boston, 2007.

\bibitem[Brooks(2007)]{Brooks2007}
William Brooks.
\newblock Pragmatics of silence.
\newblock In Nicky Losseff and Jenny Doctor, editors, \emph{Silence, Music,
  Silent Music}, pages 97--126. Ashgate, Aldershot, UK, 2007.

\bibitem[Cage and Kuhn(2016)]{Cage2016}
J.~Cage and L.~Kuhn.
\newblock \emph{The Selected Letters of {J}ohn {C}age}.
\newblock Wesleyan University Press, Middletown, CT, 2016.

\bibitem[Cage(1961)]{Cage1961}
John Cage.
\newblock \emph{Silence}.
\newblock Wesleyan University Press, Middletown, CT, 1961.

\bibitem[Cage(1962)]{Cage1962}
John Cage.
\newblock \emph{0'00'' (4'33'' No. 2)}.
\newblock Henmar Press, New York, 1962.

\bibitem[Cage(1968)]{Cage1968}
John Cage.
\newblock \emph{A Year from Monday: Lectures and Writings}.
\newblock Calder and Boyars, New York, 1968.

\bibitem[Cage(1986)]{Cage1986}
John Cage.
\newblock \emph{4'33''}.
\newblock Peters Edition, New York, 1986.

\bibitem[Cage(1989)]{Cage1989}
John Cage.
\newblock \emph{$\textrm{One}^{3}$}.
\newblock John Cage Trust, New York, 1989.
\newblock unpublished score, Fax manuscript.

\bibitem[Craenen(2014)]{Craenen2014}
Paul Craenen.
\newblock \emph{Composing under the Skin: The Music-making Body at the
  Composer's Desk}.
\newblock Leuven University Press, Leuven, 2014.

\bibitem[Daniels(2016)]{Daniels2016}
Dieter Daniels.
\newblock Silence and the void: Aesthetics of absence in space and time.
\newblock In Y.~Kaduri, editor, \emph{The Oxford Handbook of Sound and Image in
  Western Art}, Oxford Handbooks. Oxford University Press, Oxford, 2016.
\newblock \url{https://books.google.de/books?id=29f\_DQAAQBAJ}.

\bibitem[Davies(1997)]{Davies1997}
Stephen Davies.
\newblock John {C}age's 4'33": Is it music?
\newblock \emph{Australasian Journal of Philosophy}, 75\penalty0 (4):\penalty0
  448--482, 1997.

\bibitem[Davies(2003)]{Davies2003}
Stephen Davies.
\newblock \emph{Themes in the Philosophy of Music}.
\newblock Oxford University Press, Oxford, 2003.

\bibitem[Dodd(2018)]{Dodd2018}
Julian Dodd.
\newblock What 4'33'' is.
\newblock \emph{Australasian Journal of Philosophy}, 96\penalty0 (4):\penalty0
  629--641, 2018.

\bibitem[Fetterman(1996)]{Fetterman1996}
William Fetterman.
\newblock \emph{John Cage's theatre pieces: notations and performances}.
\newblock Harwood Academic Publishers GmbH, Amsterdam, 1996.

\bibitem[Fowler(2012)]{Fowler2012}
Michael Fowler.
\newblock On the relationship between inside and outside: Conceptualising
  acoustic space in {J}ohn {C}age's {V}ariations {IV}.
\newblock \emph{Architectural Theory Review}, 17\penalty0 (1):\penalty0
  158--177, 2012.

\bibitem[Fowler(2019)]{Fowler2019}
Michael Fowler.
\newblock John {C}age's {S}ilent {P}iece and the {J}apanese gardening technique
  of \textit{shakkei}: Formalizing {W}hittington's conjecture through
  conceptual graphs.
\newblock \emph{Journal of Mathematics and Music}, 13\penalty0 (1):\penalty0
  4--26, 2019.

\bibitem[Foy(2010)]{Foy2010}
George~Michelsen Foy.
\newblock \emph{Zero Decibels: The Quest for Absolute Silence}.
\newblock Simon and Schuster, New York, 2010.

\bibitem[Gann(2010)]{Gann2010}
Kyle Gann.
\newblock \emph{No Such Thing as Silence: {J}ohn {C}age's \emph{4'33''}}.
\newblock Yale University Press, New Haven, CT, 2010.

\bibitem[Gentry(2018)]{Gentry2018}
Philip Gentry.
\newblock The cultural politics of 4'33'': identity and sexuality.
\newblock In Matthieu Saladin, editor, \emph{Tacet \#1: Who is John Cage?} Les
  presses du r{\'e}el, Dijon, 2018.

\bibitem[Greimas(1983)]{Greimas1983}
A.~J. Greimas.
\newblock \emph{Structural Semantics: An Attempt at a Method}.
\newblock University of Nebraska Press, Lincoln, 1983.

\bibitem[Johnson(1970)]{Johnson1970}
J.~Johnson.
\newblock There is no silence now.
\newblock In J.~Cage and R.~Kostelanetz, editors, \emph{John {C}age: An
  anthology}, pages 145--149. Da Capo Press, New York, 1970.

\bibitem[Joseph(1997)]{Joseph1997}
Branden~W. Joseph.
\newblock John {C}age and the architecture of silence.
\newblock \emph{October}, 81\penalty0 (Summer):\penalty0 80--104, 1997.

\bibitem[Kahn(1997)]{Kahn1997}
Douglas Kahn.
\newblock John {C}age: Silence and silencing.
\newblock \emph{The Musical Quarterly}, 81\penalty0 (4):\penalty0 556--598,
  1997.

\bibitem[Katschthaler(2016)]{Katschthaler2016}
Karl Katschthaler.
\newblock Absence, presence and potentiality: {J}ohn {C}age's \emph{4'33''}
  revisited.
\newblock In Werner Wolf and Walter Bernhart, editors, \emph{Silence and
  Absence in Literature and Music}, pages 166--179. Rodopi, Boston, 2016.

\bibitem[Kim-Cohen(2009)]{Kim-Cohen2009}
Seth Kim-Cohen.
\newblock \emph{In the blink of an ear: towards a non-cochlear sonic art}.
\newblock Contimuum, New York, 2009.

\bibitem[Kostelanetz and Cage(2003)]{KostelanetzCage2003}
R.~Kostelanetz and J.~Cage.
\newblock \emph{Conversing with Cage}.
\newblock Readers' Guides to Essential Criticism. Routledge, London, 2003.

\bibitem[Kotz(2009)]{Kotz2009}
Liz Kotz.
\newblock Cagean structures.
\newblock In Julia Robinson, editor, \emph{The Anarchy of Silence: John Cage
  and experimental art}. Musea d'{Art} {C}ontemporani de {B}arcelona,
  Barcelona, 2009.

\bibitem[Kruse(2019)]{Kruse2019}
Robert Kruse.
\newblock John {C}age and nonrepresentational spaces of music.
\newblock \emph{Social \& Cultural Geography}, 2019.
\newblock \url{10.1080/14649365.2019.1628291}.

\bibitem[LaBelle(2015)]{Labelle2015}
Brandon LaBelle.
\newblock \emph{Background Noise, Second Edition: Perspectives on Sound Art}.
\newblock Bloomsbury Publishing, London, 2015.

\bibitem[MacLane(2013)]{Maclane2013}
Saunders MacLane.
\newblock \emph{Categories for the Working Mathematician}.
\newblock Graduate Texts in Mathematics. Springer, New York, 2013.

\bibitem[Massey(2009)]{Massey2005}
D.~Massey.
\newblock \emph{For space}.
\newblock Sage, London, 2009.

\bibitem[Moholy-Nagy and Hoffmann(2012)]{Moholy-Nagy2012}
L.~Moholy-Nagy and D.M. Hoffmann.
\newblock \emph{The New Vision: Fundamentals of Bauhaus Design, Painting,
  Sculpture, and Architecture}.
\newblock Dover Fine Art, History of Art. Dover Publications, London, 2012.

\bibitem[Pisaro(2018)]{Pisaro2018}
Michael Pisaro.
\newblock Nicht alles, nicht nichts, etwas. {T}he opposing tensions of {J}ohn
  {C}age's 0'00'' and {R}oaratorio.
\newblock In Matthieu Saladin, editor, \emph{Tact:\#1: Who is John Cage?},
  pages 110--125. Les presses du r\'{e}el, Dijon, 2018.

\bibitem[Pritchett(1996)]{Pritchett1996}
James Pritchett.
\newblock \emph{The Music of {J}ohn {C}age}.
\newblock Cambridge University Press, Cambridge, 1996.

\bibitem[Pritchett(2018)]{Pritchett2018}
James Pritchett.
\newblock Silence changed: {O}ne$^{3}$, 2018.
\newblock
  \url{http://rosewhitemusic.com/piano/2018/10/01/silence-changed-one3/}.

\bibitem[Retallack(2013)]{Retallack2013}
Joan Retallack.
\newblock Poethics of a complex realism.
\newblock In Marjorie Perloff and Charles Junkerman, editors, \emph{John Cage:
  Made in America}, pages 242--274. University of Chicago Press, Chicago, 2013.

\bibitem[Schubert(1972)]{Schubert1972}
Horst Schubert.
\newblock \emph{Categories}.
\newblock Springer-Verlag, New York, 1972.
\newblock Eva Gray, translator.

\bibitem[Shultis(2013)]{Shultis2013}
Christopher Shultis.
\newblock \emph{Silencing the Sounded Self: {J}ohn {C}age and the {A}merican
  {E}xperimental {T}radition}.
\newblock University Press of New England, Lebanon, NH, 2013.

\bibitem[Spivak(2010)]{Spivak2013b}
David~I. Spivak.
\newblock Functorial data migration.
\newblock \emph{CoRR}, abs/1009.1166, 2010.
\newblock \url{http://arxiv.org/abs/1009.1166}.

\bibitem[Spivak(2013{\natexlab{a}})]{Spivak2013}
David~I Spivak.
\newblock Category theory for scientists (old version).
\newblock \emph{\url{arXiv preprint arXiv:1302.6946}}, 2013{\natexlab{a}}.

\bibitem[Spivak(2013{\natexlab{b}})]{Spivak2013c}
David~I. Spivak.
\newblock Database queries and constraints via lifting problems.
\newblock \emph{Mathematical Structures in Computer Science}, 24\penalty0 (6),
  2013{\natexlab{b}}.
\newblock ISSN 1469-8072.
\newblock \url{http://dx.doi.org/10.1017/S0960129513000479}.

\bibitem[Spivak and Kent(2012)]{Spivak2012}
David~I Spivak and Robert~E Kent.
\newblock Ologs: a categorical framework for knowledge representation.
\newblock \emph{PLoS One}, 7\penalty0 (1):\penalty0 e24274, 2012.

\bibitem[Swed(1993)]{Swed1993}
Mark Swed.
\newblock John {C}age: {S}eptember 5, 1912-{A}ugust 12, 1992.
\newblock \emph{The Musical Quarterly}, 77\penalty0 (1):\penalty0 132--144,
  1993.

\bibitem[Vandendorpe(1993)]{Vandendorpe1993}
Christian Vandendorpe.
\newblock Actant.
\newblock In I.R. Makaryk, editor, \emph{Encyclopedia of Contemporary Literary
  Theory: Approaches, Scholars, Terms}, Theory / culture, page 505. University
  of Toronto Press, Toronto, 1993.

\bibitem[Whittington(2013)]{Whittington2013}
Stephen Whittington.
\newblock Digging in {J}ohn {C}age's garden: {C}age and {R}y\={o}anji.
\newblock \emph{Malaysian Music Journal}, 2\penalty0 (2):\penalty0 12--21,
  2013.

\end{thebibliography}
\bibliographystyle{plainnat}
\end{document}